\documentclass[preprint,11pt]{elsarticle}

\usepackage[utf8]{inputenc}
\usepackage[english]{babel}
\usepackage[hmargin=0.9in]{geometry}
\usepackage{amsmath}
\usepackage{amsfonts}
\usepackage{amssymb}
\usepackage{graphicx}
\usepackage[colorlinks=true,linkcolor=blue]{hyperref}
\usepackage{bm}
\usepackage{bbold}
\usepackage{array}
\usepackage{subcaption}
\usepackage{siunitx}
\usepackage{bbm}
\usepackage[capitalise,noabbrev,nameinlink]{cleveref}
\usepackage{appendix}
\usepackage{booktabs}
\usepackage{placeins} 

\date{March 12, 2026}

\newcommand{\D}{\mathrm{d}}
\newcommand{\pd}[2]{\frac{\partial #1}{\partial #2}}
\newcommand{\dx}{\,\mathrm{d}\mathbf{x}}
\newcommand{\td}[2]{\frac{\D #1}{\D #2}}
\newcommand{\hv}{\left(h\mathbf{v}\right)}
\newcommand{\bmat}{M_{\Gamma}}
\newcommand{\mmin}{\min\phantom{a}\hspace{-1.8mm}}
\newcommand{\mmax}{\max\phantom{i}\hspace{-1.8mm}}

\newcommand{\tvd}{(TVD)}
\newcommand{\Lonelone}{\texorpdfstring{\ensuremath{L^1 \left( l_1 \right)}}{L1(l1)}}
\newcommand{\Loneltwo}{\texorpdfstring{\ensuremath{L^1 \left( l_2 \right)}}{L1(l2)}}

\DeclareMathOperator*{\minimize}{minimize}
\DeclareMathOperator*{\maximize}{maximize}

\makeatletter
\def\ps@pprintTitle{%
  \let\@oddhead\@empty
  \let\@evenhead\@empty
  \def\@oddfoot{
    \footnotesize\itshape
    \hfill\@date
  }%
  \let\@evenfoot\@oddfoot}
\makeatother

\begin{document}

\begin{frontmatter}
    \title{Bathymetry reconstruction via optimal control in well-balanced finite element methods for the shallow water equations}
    \author{Falko Ruppenthal\texorpdfstring{\corref{cor}}{}}
    \ead{falko.ruppenthal@math.tu-dortmund.de}
    \cortext[cor]{Corresponding author}
    \author{Dmitri Kuzmin}
    \ead{kuzmin@math.tu-dortmund.de}
    \address{Institute of Applied Mathematics (LS III), TU Dortmund University, Vogelpothsweg 87, D-44227 Dortmund, Germany}

    \begin{abstract}
        Accurate prediction of shallow water flows relies on precise bottom topography data, yet direct bathymetric surveys are expensive and time-consuming. In contrast, remote sensing platforms such as radar or satellite altimetry provide accurate free surface observations. This disparity motivates a data-driven reconstruction strategy: invert the shallow water equations to estimate the bathymetry that yields the best fit to the governing dynamics. We introduce a new direct reconstruction technique that extracts bathymetric features from widely available free surface measurements. The underlying inverse problem of determining an unknown bathymetry profile from observed wave elevations is inherently ill-posed. Small perturbations in the data may lead to large deviations in the reconstructed topography, and discontinuities or sharp gradients further exacerbate instability. To stabilize the inversion, we formulate an optimal-control problem, wherein a cost functional penalizes deviations between simulated and measured free surface elevation while enforcing a state equation for the flow dynamics. To suppress noise and preserve sharp depth variations, the framework is augmented with $L^1$ regularization and total variation denoising. These sparsity-promoting terms encourage piecewise-smooth solutions, allowing changes in the bathymetry to be captured without excessive smoothing. Numerical experiments on synthetic noisy data and discontinuous bathymetry demonstrate robust performance in reconstructing unknown bathymetry.
    \end{abstract}


\end{frontmatter}

\section{Introduction}\noindent
Accurate bottom topography (bathymetry) data is fundamental to various applications, including coastal engineering, marine navigation, and climate change research. Direct sonar surveys remain cost-prohibitive and limited to localized studies, making remote sensing the preferred method for large-scale investigations of water bodies (see Smith and Sandwell \citep{smith2004conventional}). Satellite altimetry and optical imagery provide precise free surface elevation measurements. In contrast, obtaining accurate bathymetry information remains challenging. A promising alternative is to recover bottom topography from surface observations by solving an inverse problem. This approach requires carefully designed regularization techniques to stabilize the ill-conditioned inverse problem and ensure physically consistent bathymetric estimates. We present a novel framework for reconstructing bathymetry with high accuracy from synthetic or measured free surface data.

Recent literature provides a collection of reconstruction tools based on data assimilation, machine learning, adjoint-based PDE optimization, and variational methods for solving inverse problems. Durand et al. \citep{durand2008estimation} developed an ensemble-based data-assimilation method that reconstructs seafloor depth and slope from synthetically generated data. It achieves accurate reconstructions even in the presence of measurement noise, although uncertainty quantification remains a challenge.
Machine-learning strategies have also emerged in recent years. Liu et al. \citep{liu2024bathymetry} trained a deep neural network using solutions of the shallow water equations with known bathymetry to recover an unknown bathymetric profile. Their formulation successfully captures bathymetric features by incorporating both value and slope regularizations into the loss function.
Adjoint-based PDE optimization has gained popularity because it enforces conservation laws and boundary conditions explicitly. Angel et al. \citep{angel2024bathymetry} applied such an adjoint-based PDE optimization to the results of laboratory experiments conducted on an artificial beach. The bathymetry was reconstructed by minimizing a functional that quantifies the misfit between measured and simulated surface elevations. The involved adjoint problem was derived analytically from the governing equations. The approach of \citep{angel2024bathymetry} demonstrates the utility of laboratory data when the forward model is fully specified.

Montecinos et al. \citep{montecinos2019numerical} developed a coupled finite-volume solver capable of handling conservative and non-conservative hyperbolic equations. Their formulation treats forward and adjoint problems in a unified manner and reconstructs time-independent bathymetry with a gradient descent method. Lecaros et al. \citep{lecaros2024optimal} extended this framework to time-dependent bathymetry reconstruction. They established well-posedness of the resulting problem and demonstrated that their approach delivers superior accuracy when benchmarked against reference solutions.
In addition, several variational approaches have emerged. Gessese et al. \citep{gessese2011reconstruction,gessese2012direct} used synthetic data generated by a forward solver and introduced explicit one-shot methods for finite volume and finite difference schemes that are similar for the forward and inverse model, enabling fast bathymetry updates. Hajduk et al. \citep{hajduk2020bathymetry} employed a discontinuous Galerkin-continuous Galerkin (DG-CG) discretization for a bathymetry evolution equation, adding an artificial viscosity term to regularize the discrete inverse problem and mitigate spurious oscillations due to noisy measurements. Furthermore, Britton et al. \citep{britton2021recovery} framed bathymetry as a linear combination of a priori known functions and solved for a time-dependent coefficient using $L^1$ and $H^1$ regularization simultaneously. Recent work by Lamsahel and Rosier \citep{lamsahel2025direct} introduced a novel inversion strategy for channel topography by using free surface observations and their first two time derivatives at a single time instant. The algorithm is applicable to the one-dimensional shallow water equations in subcritical and supercritical flow regimes.

Another major challenge is the treatment of measurement noise. Existing strategies include total variation denoising (TVD) for data processing \citep{vogel1996iterative,kunisch2004total} and ensemble-based filtering procedures for data assimilation \citep{durand2008estimation}. Our approach employs $L^1$ optimization to promote sparsity robustly in noisy measurements. Unlike earlier works that apply a two-stage filtering-inversion pipeline, we embed an $L^1$-penalty term into the cost functional to combine regularization and optimization.

Building on the variational framework of \citep{hajduk2020bathymetry}, we formulate an optimal control problem that recovers bathymetry from free surface elevation. Specifically, we use the flux-potential methodology of Ruppenthal et al. \citep{ruppenthal2023optimal,ruppenthal2024scalable} for bathymetry reconstruction purposes, taking advantage of its inherent locality and conservation properties. We derive the first-order optimality conditions and outline several numerical solution strategies. We further enrich the optimization-based framework with convex (TVD) and non-smooth ($L^1$) regularizations.

This paper proceeds as follows. \Cref{sec:governing} presents the shallow water equations and establishes the mathematical setting of the problem. \Cref{sec:numerical} gives an overview of the low-order algebraic Lax-Friedrichs scheme for the forward problem and of the high-order extension using monolithic convex limiting. Subsequently, \Cref{sec:inverse_opt} formulates the bathymetry reconstruction as a coupled inverse problem, introduces the required notation, discusses well-posedness issues arising from the underlying hyperbolic system, and details the discrete state equation formulation for optimization (including the cost functional and optimality conditions). In \Cref{sec:regularization}, we describe oscillation-suppression strategies based on total variation denoising and $L^1$ regularization that are incorporated into the formulation of optimization problems. In \Cref{sec:implementation}, we detail the numerical solution strategy including convergence acceleration and software implementation. Furthermore, \Cref{sec:experiments} reports numerical experiments on one- and two-dimensional benchmarks, illustrating the robustness and accuracy of the proposed method and providing an analysis of its sensitivity to measurement noise.

\section{Governing equations}\label{sec:governing}\noindent
In this section, we formulate the shallow water equations and outline the key features of the associated inverse problem. Let $\Omega \subset \mathbb{R}^d \, \left(d = 1,2\right)$ denote a bounded domain with boundary $\Gamma = \partial\Omega$. The Saint-Venant (shallow water) equations with variable bathymetry can be written as
\begin{subequations}\label{SWE:both}
    \begin{align}
        \frac{\partial h}{\partial t} + \nabla\! \cdot\, \hv &= 0, &\text{in }\Omega\times I, \label{SWE:CON} \\
        \frac{\partial \hv}{\partial t} + \nabla\! \cdot\! \left( \frac{\hv \otimes \hv}{h} + \frac12\, g\,h^2\, \mathbf{I}_d \right) + g\, h\, \nabla \tilde b &= 0, &\text{in }\Omega\times I. \label{SWE:MOM}
    \end{align}
\end{subequations}
The conserved variables are the water height $h: \,\Omega\times I \rightarrow \mathbb R$ and the depth-averaged momentum $\hv: \, \Omega\times I \rightarrow \mathbb R^d$. The vector field $\mathbf v=\hv/h$ corresponds to the horizontal velocity. The bathymetry is assumed to be constant in time. Therefore, it is represented by a time-independent function $\tilde b:\,\Omega\rightarrow \mathbb R$. In this context, $g$ denotes gravitational acceleration, and $\mathbf{I}_d\in\mathbb{R}^{d \times d}$ represents the identity tensor. The time variable is denoted by $t$. We solve the nonlinear hyperbolic system \eqref{SWE:both} over the time interval $I=\left[0,T\right]$, where $T > 0$. The initial conditions and external data of the weakly imposed boundary conditions are provided by the functions $h_e$ and $\hv_e$.

\section{Forward numerical methodology}\label{sec:numerical}\noindent
A critical prerequisite for any inverse problem framework is the selection of a robust forward solver. For a system of hyperbolic balance laws, such a solver must deliver high-order accuracy while remaining stable under realistic wave propagation scenarios. Recent advances in this area have been documented extensively \citep{audusse2004fast,kurganov2007second,audusse2016kinetic,guermond2018well,hajduk2022algebraically}.
Our framework requires a forward scheme capable of accommodating bathymetric effects through an adjustable source term. The employed finite element discretization of the shallow water equations is based on the methodology developed in \citep{hajduk2022algebraically,hajduk2022bound}. It ensures three essential properties: positivity preservation for the water height, well-balancedness (preservation of lakes at rest), and entropy stability. These properties are enforced using a tailored monolithic convex limiting (MCL) algorithm. For a comprehensive discussion on modern limiting strategies for hyperbolic conservation laws, we refer the reader to the recent review by Kuzmin and Hajduk \citep{kuzmin2023property}.

\subsection{Finite-element discretization}\noindent
Let $\varphi_i: \,\overline{\Omega} \rightarrow \left[0,1\right]$ denote the nodal Lagrange basis function associated with node $i$. A continuous finite element approximation $u_h\in V_h \subset H^1(\Omega)$ to an exact weak solution $u=(h,\hv)$ of the hyperbolic system \eqref{SWE:both} is defined by
\begin{align*}
    u_h (\mathbf x,t) = \sum_{i=1}^{N_h} u_i(t) \varphi_i (\mathbf x).
\end{align*}
We denote by
\begin{align*}
    m_{ij} = \int_{\Omega} \varphi_i \varphi_j \dx, \quad m_i = \sum^{N_h}_{j=1} m_{ij}
\end{align*}
the entries of the consistent mass matrix $M_C = \left(m_{ij}\right)_{i,j=1}^{N_h}$ and the diagonal entries of the lumped matrix $M_L = \text{diag}\left(m_1,\dots,m_{N_h}\right)$, which can be interpreted as a quadrature-based approximation to $M_C$.
The coefficients of the discrete gradient operator $\mathbf C=\left(\mathbf c_{ij}\right)_{i,j=1}^{N_h}$ are given by
\begin{align*}
    \mathbf{c}_{ij} = \int_{\Omega} \varphi_i \nabla \varphi_j \dx.
\end{align*}
We store the indices of basis functions $\varphi_j$ that have overlapping support with $\varphi_i$ in $\mathcal N_i = \left\{ j: m_{ij} \neq 0 \right\}$ and define $\mathcal N^*_i = \mathcal N_i \backslash \left\{i\right\}$.

Temporal discretization employs a second-order explicit strong stability-preserving Runge-Kutta (Heun) scheme. Under appropriate CFL conditions, this scheme preserves the invariant domain of the shallow water system.
We impose the boundary conditions with a Rusanov flux that uses the internal states of the conserved variables and the prescribed external data ($h_e, \hv_e$). The numerical viscosity coefficient is estimated via the maximum eigenvalue of the shallow water equations.

We next review the numerical strategies designed in \citep{hajduk2022algebraically,hajduk2022bound} to discretize and solve the shallow water equations. We begin by presenting a low-order scheme that approximates the local Lax-Friedrichs flux. Next, we outline a high-order approach based on monolithic convex limiting. Finally, we analyze the inverse problem of reconstructing bathymetry from synthetic data.

\subsection{Low-order (ALF) scheme}\noindent
For the forward problem, we discretize the Saint-Venant system with an algebraic Lax-Friedrichs (ALF) scheme that treats the bathymetry source term explicitly in time. For homogeneous hyperbolic systems, this local Lax-Friedrichs scheme was analyzed in \cite{guermond2016invariant}. It was subsequently modified for shallow water dynamics via hydrostatic reconstruction in \cite{azerad2017well}. An alternative way to preserve lake-at-rest configurations and enforce positivity preservation through limiting techniques was proposed in \cite{hajduk2022algebraically,hajduk2022bound}. Stability of these extensions hinges on the artificial diffusion coefficients $d_{ij}$. They are defined for every edge ($i,j$) of the mesh using an upper bound for the maximum wave speed of the one-dimensional Riemann problem corresponding to the normal direction $\mathbf n_{ij}=\mathbf c_{ij}/|\mathbf c_{ij}|$. For the shallow water equations, the relevant wave speeds are $\lambda_1 = \mathbf{v} \cdot \mathbf{n}$ and $\lambda_{2,3} = \mathbf{v}\cdot \mathbf{n} \pm \sqrt{gh}$. The corresponding upper bound $|\mathbf{v}\cdot \mathbf{n}|+\sqrt{gh}$ is used in the definition
\begin{align*}
  d_{ij} = \max&\left\{
  |\mathbf{v}_i \cdot \mathbf{c}_{ij}| + |\mathbf{c}_{ij}|\sqrt{gh_i},
  |\mathbf{v}_i \cdot \mathbf{c}_{ji}| + |\mathbf{c}_{ji}|\sqrt{gh_i},\right. \\ &\hphantom{\big\{}\left.
  |\mathbf{v}_j \cdot \mathbf{c}_{ij}| + |\mathbf{c}_{ij}|\sqrt{gh_j},
  |\mathbf{v}_j \cdot \mathbf{c}_{ji}| + |\mathbf{c}_{ji}|\sqrt{gh_j} \right\}
\end{align*}
of the artificial diffusion coefficient associated with the edge
($i,j$).

The low-order ALF method approximates system \eqref{SWE:both} by \citep{hajduk2022algebraically,hajduk2022bound}
\begin{align}
    m_i \td{h_i}{t} &= \sum_{j \in\mathcal N^*_i} \left[ d_{ij} \left( h_j - h_i + \alpha^b_{ij} \left( b_j - b_i \right) \right) - \left(\hv_j - \hv_i \right) \cdot \mathbf{c}_{ij} \right], \label{eq:alf_con}\\
    m_i \td{\!\hv_i}{t} &= \sum_{j \in\mathcal N^*_i} \left[ d_{ij} \left( \hv_j - \hv_i + \frac{\alpha^b_{ij}}{2} \left(b_j - b_i \right) \left( \mathbf{v}_i + \mathbf{v}_j \right) \right) \right. \nonumber\\
    &\hphantom{= \sum_{j \in\mathcal N^*_i} \Big [}\left. - \left(\mathbf{f}_j - \mathbf{f}_i \right) \mathbf{c}_{ij} - \frac{\alpha^b_{ij}}{2} g \left(h_j + h_i\right) \left(b_j - b_i \right) \mathbf{c}_{ij} \right].\nonumber
\end{align}
The nodal flux $\mathbf f_i=\mathbf f(u_i)$ appearing in the momentum balance is defined as
\begin{align*}
    \mathbf{f}_i = \frac{\hv_i \otimes \hv_i}{h_i} + \frac12\, g\,h^2_i\, \mathbf{I}_d.
\end{align*}
The correction factor $\alpha^b_{ij} \in [0,1]$ was introduced in \cite{hajduk2022algebraically} to guarantee positivity preservation for $h$. Given the focus on deep-water bathymetry reconstruction, $\alpha^b_{ij}$ is typically set to unity in our work.

\subsection{High-order (MCL) scheme}\noindent
The local Lax-Friedrichs scheme exhibits significant numerical diffusion, so we seek to recover second-order accuracy in smooth regions while rigorously enforcing local discrete maximum principles and positivity preservation. Kuzmin's MCL strategy \cite{kuzmin2020monolithic} provides a systematic framework for enforcing local and global bounds in the context of hyperbolic (systems of) conservation laws. Hajduk \cite{hajduk2022algebraically} adapted the MCL approach to the shallow water equations with non-flat bathymetry. The central idea is to generalize the concept of flux-corrected bar states to systems with source terms. The ALF evolution equation for $h$ may be reformulated as
\begin{align} \label{eq:llf}
    m_i \td{h_i}{t} = \sum_{j \in \mathcal N^*_i} 2 d_{ij} \left( \overline{h}^b_{ij} - h_i \right),
\end{align}
where
\begin{align*}
    \overline{h}^b_{ij} = \frac{h_i + h_j}{2} - \frac{\left(\hv_j - \hv_i \right) \cdot \mathbf{c}_{ij}}{2 d_{ij}} + \frac{\alpha^b_{ij}}{2} \left( b_j - b_i \right).
\end{align*}
To recover second-order accuracy, the antidiffusive fluxes
\begin{align*}
    f^h_{ij} = m_{ij} \left( \dot h_i - \dot h_j \right) + d_{ij} \left( h_i - h_j + \alpha^b_{ij} \left( b_i - b_j \right) \right)
\end{align*}
are added on the right-hand side of the target scheme
\begin{align*}
    m_i \td{h_i}{t} = \sum_{j \in \mathcal N^*_i} \left[ 2 d_{ij} \left( \overline{h}^b_{ij} - h_i \right) + f^h_{ij} \right].
\end{align*}
The first term in the formula for $f^h_{ij}$ eliminates the mass lumping error, while the second term offsets the artificial diffusion inherent in the low-order scheme. Our goal is to remove as much of this diffusion as possible without producing physically or numerically unacceptable states.

The sequential MCL scheme presented in \citep{hajduk2022algebraically} enforces the inequality constraints
\begin{align*}
     h^{\min}_i \leq \overline{h}^{b*}_{ij} = \overline{h}^b_{ij} + \frac{f^{h*}_{ij}}{2 d_{ij}} \leq h^{\max}_i.
\end{align*}
For this particular problem, a closed-form solution exists and can be written as
\begin{equation*}
    f^{h*}_{ij} = \begin{cases}
        \begin{aligned}
        \mmin &\left\{f^h_{ij}, 2 d_{ij} \mmin \left\{ h^{\max}_i - \overline{h}^b_{ij}, \overline{h}^b_{ji} - h^{\min}_j \right\} \right\} \quad&\text{if}\; f^h_{ij} \geq 0, \\
        \mmax &\left\{f^h_{ij}, 2 d_{ij} \mmax \left\{ h^{\min}_i - \overline{h}^b_{ij}, \overline{h}^b_{ji} - h^{\max}_j \right\} \right\} \quad&\text{if}\; f^h_{ij} < 0.
        \end{aligned}
    \end{cases}
\end{equation*}
The MCL constraints for $\overline{h}^{b*}_{ij}$ are feasible if $h^{\min}_i \leq \overline{h}^b_{ij}\leq h^{\max}_i$. Hence, the data set for constructing the upper and lower bounds should include the low-order bar states \cite{hajduk2022algebraically}. To make the bounds
\begin{align*}
    h^{\max}_i = \max\left\{ \max_{j \in \mathcal N_i} h_j, \max_{j \in \mathcal N^*_i} \overline{h}^b_{ij}\right\} \quad\text{and}\quad h^{\min}_i = \min\left\{ \min_{j \in \mathcal N_i} h_j, \min_{j \in \mathcal N^*_i} \overline{h}^b_{ij}\right\}
\end{align*}
less restrictive, we additionally include the nodal states with indices $j \in \mathcal N_i$. The flux-corrected bar states $h^{b*}_{ij}$ replace the low-order ones $h^b_{ij}$ in \eqref{eq:llf} to construct a limited high-order approximation.

After limiting the water height in the above manner, we express the target discretization of the momentum equation in terms of the bar states
\begin{align*}
    \overline{\hv}^b_{ij} = \frac{\hv_i + \hv_j}{2} - \frac{\left(\mathbf{f}_j - \mathbf{f}_i + \frac{\alpha^b_{ij}}{2} g \left( h_i + h_j \right) \left( b_j - b_i \right) \mathbf{I}_d \right) \mathbf{c}_{ij}}{2 d_{ij}}
    + \frac{\alpha^b_{ij}}{4} \left(b_j - b_i \right) \left( \mathbf{v}_i + \mathbf{v}_j \right)
\end{align*}
and the antidiffusive fluxes
\begin{align*}
    \mathbf{f}^{\hv}_{ij} = m_{ij} \left( \dot{\hv}_i - \dot{\hv}_j \right) + d_{ij} \left( \hv_i - \hv_j + \frac{\alpha^b_{ij}}{2} \left( b_i - b_j \right) \left(\mathbf{v}_i + \mathbf{v}_j \right) \right).
\end{align*}
The high-order target
\begin{align*}
    m_i \td{\!\hv_i}{t} = \sum_{j \in \mathcal N^*_i} \left[ 2 d_{ij} \left( \overline{\hv}^b_{ij} - \hv_i \right) + \mathbf{f}^{\hv}_{ij} \right]
\end{align*}
has the same form as the semi-discrete equation for the water height.

To formulate MCL constraints for the velocity $\mathbf v$, we define the states (cf. \citep{dobrev2018sequential,kuzmin2020monolithic,hajduk2019new})
\begin{align*}
    \overline{h}^*_{ij} = \overline{h}^{b*}_{ij} - \frac{\alpha^b_{ij}}{2} \left( b_j - b_i \right),\qquad \overline{\mathbf{v}}_{ij} = \frac{\overline{\hv}^b_{ij} + \overline{\hv}^b_{ji}}{\overline{h}^b_{ij} + \overline{h}^b_{ji}}
\end{align*}
and the auxiliary fluxes
\begin{align*}
    \mathbf{g}^{\hv}_{ij} = \mathbf{f}^{\hv}_{ij} + 2 d_{ij} \left( \overline{\hv}^b_{ij} - \overline{h}^*_{ij} \overline{\mathbf{v}}_{ij} \right).
\end{align*}
The constraints imposed on the momentum bar states are given by
\begin{align*}
    \overline{h}^*_{ij} \mathbf{v}^{\min}_i \leq \overline{\hv}^{b*}_{ij} = \overline{\hv}^b_{ij} + \frac{\mathbf{f}^{\hv *}_{ij}}{2 d_{ij}} &= \overline{h}^*_{ij} \overline{\mathbf{v}}_{ij} + \frac{\mathbf{g}^{\hv *}_{ij}}{2 d_{ij}} \leq \overline{h}^*_{ij} \mathbf{v}^{\max}_i.
\end{align*}
Here and below, inequalities involving vectors are meant to hold for individual components. The upper and lower bounds
\begin{align*}
    \mathbf{v}^{\max}_i &= \mmax\left\{ \underset{j \in \mathcal N_i}{\mmax}\; \mathbf{v}_j, \underset{j \in \mathcal N^*_i}{\mmax}\; \overline{\mathbf{v}}_{ij}, \underset{j \in \mathcal N^*_i}{\mmax}\; \frac{\overline{\hv}^b_{ij}}{\overline{h}_{ij}} \right\},\\
    \mathbf{v}^{\min}_i &= \mmin\left\{ \underset{j \in \mathcal N_i}{\mmin}\; \mathbf{v}_j, \underset{j \in \mathcal N^*_i}{\mmin}\; \overline{\mathbf{v}}_{ij}, \underset{j \in \mathcal N^*_i}{\mmin}\; \frac{\overline{\hv}^b_{ij}}{\overline{h}_{ij}} \right\}
\end{align*}
for the velocity components are enforced using
\begin{align*}
    \mathbf{g}^{\hv*}_{ij} = \begin{cases} \begin{aligned}
    \mmin \left\{\mathbf{g}^{\hv}_{ij}, 2 d_{ij} \mmin \left\{ \overline{h}^*_{ij} \left( {\mathbf{v}}^{\mmax}_i - \overline{\mathbf{v}}_{ij} \right), \overline{h}^*_{ji} \left( {\overline{\mathbf{v}}_{ij} - \mathbf{v}}^{\mmin}_j \right) \right\} \right\} \quad&\text{if}\; \mathbf{g}^{\hv}_{ij} \geq 0, \\
    \mmax \left\{\mathbf{g}^{\hv}_{ij}, 2 d_{ij} \mmax \left\{ \overline{h}^*_{ij} \left( {\mathbf{v}}^{\mmin}_i - \overline{\mathbf{v}}_{ij} \right), \overline{h}^*_{ji} \left( {\overline{\mathbf{v}}_{ij} - \mathbf{v}}^{\mmax}_j \right) \right\} \right\} \quad&\text{if}\; \mathbf{g}^{\hv}_{ij} < 0
    \end{aligned}\end{cases}
\end{align*}
to construct 
\begin{align*}
    \mathbf{f}^{\hv*}_{ij} = \mathbf{g}^{\hv*}_{ij} - 2 d_{ij} \left( \overline{\hv}^b_{ij} - \overline{h}^*_{ij} \overline{\mathbf{v}}_{ij} \right).
\end{align*}

In summary, the ALF scheme ensures nonlinear stability and preservation of invariant domains by introducing strong artificial diffusion. The MCL version attempts to recover the baseline Galerkin discretization by adding bound-preserving antidiffusive corrections. This adjustment yields second-order accuracy in smooth regions while still enforcing positivity of the water height. Entropy stability can also be ensured using a limiter-based fix \cite{hajduk2022algebraically,hajduk2022bound}. Consequently, ALF is robust but overly diffusive, whereas MCL delivers high resolution without sacrificing any essential physical properties.

\section{Inverse problem formulation and optimization}\label{sec:inverse_opt}\noindent
Building upon the robust forward solver established in \Cref{sec:numerical}, we now formulate an optimal control problem for bathymetry reconstruction. This section details the necessary modifications of the numerical scheme required for the inverse setting, specifically addressing stability issues arising from bathymetry-dependent diffusion. We derive the governing equations for the stationary bathymetry and introduce a pseudo-time stepping method for marching the solution of the ill-posed inverse problem to a steady state. Furthermore, we define the objective functional incorporating data misfit and regularization terms, and subsequently derive the first-order optimality conditions.

\subsection{Modifications to the forward solver for inverse applications}\noindent
To ensure numerical robustness within the inverse framework, we implement two specific modifications to the forward solver. First, for subcritical outflow conditions, we prescribe the external water depth as $ h^*_{\text{ext}} = h_{\text{ext}} - b_{\text{in}}$, where $h_{\text{ext}}$ denotes the reference water height outside the domain and $b_{\text{in}}$ represents the local bathymetry elevation. This adjustment mitigates spurious oscillations that arise from the interaction between internal bathymetric waves and the boundary during pseudo-time integration. Second, we address numerical artifacts generated by coupling bathymetry-dependent artificial diffusion to the water height equation (Eq.~\eqref{eq:alf_con}). The inclusion of such terms induces smearing near steep gradients and degrades convergence rates. To suppress this artifact, we exclude bathymetry-dependent diffusion from the discretized continuity equation.

\subsection{Inverse problem formulation}\noindent
Our objective is to recover a stationary bathymetry $\tilde b$ from space-time observations of the free surface elevation $H:\,\Omega\times I \rightarrow \mathbb R$. Since equation \eqref{SWE:MOM} contains only $\nabla \tilde b$, the inverse problem is ill posed: $\tilde b$ can be recovered only up to an additive constant. Imposing a prescribed value of $\tilde b$ on the inflow part of $\Gamma$ eliminates this null space.
The free surface elevation $H$ is related to water height $h$ and bathymetry $\tilde b$ by
\begin{align*}
    H = h + \tilde b \quad \text{in }\Omega\times I.
\end{align*}
Since $\tilde b$ is time-independent, this equation indicates that $\partial_t H = \partial_t h$. Substituting $\partial_t H = \partial_t h$ into equation \eqref{SWE:CON} yields the divergence constraint
\begin{align} \label{eq:NS}
    \nabla\!\cdot \hv = -\,\partial_t H.
\end{align}
System \eqref{SWE:MOM},\eqref{eq:NS} is reminiscent of the incompressible Navier-Stokes equations in a pressure-velocity formulation, wherein the bathymetry $\tilde b$ acts as a Lagrange multiplier constraining $\nabla\!\cdot\hv$.

A direct one-shot optimization requires solving the full nonlinear PDE system with respect to $b$ at each iteration of the outer loop, which significantly increases computational cost. Instead, we interpret the divergence constraint \eqref{eq:NS} as the steady state of the evolution equation
\begin{align} \label{eq:TE}
    \partial_\xi b = \partial_t H + \nabla\!\cdot \hv
\end{align}
for a tentative bathymetry $b:\,\Omega\times I \rightarrow \mathbb R$ depending on a fictitious time $\xi \ge 0$. The left-hand side of \eqref{eq:TE} vanishes as soon as $b$ converges to a steady state $\tilde b$ independent of $\xi$. All fields are discretized using a continuous, piecewise bilinear approximation on a quadrilateral mesh.

For simplicity, we identify the fictitious time variable $\xi$ with the physical simulation time $t$ in what follows. Taking the time derivative of the identity
\begin{align*}
    H = h + b,
\end{align*}
we obtain the evolution equation
\begin{align*}
    \pd{b}{t} = \pd{H}{t} - \pd{h}{t}.
\end{align*}
The stationary bathymetry $\tilde b$ is the corresponding steady state. Discretizing the above equation in space with linear finite elements and advancing $b$ in time by a single-step scheme yields
\begin{align*}
    M_C \frac{b^{n+1} - b^n}{\Delta t} = M_C \frac{H^{n+1} - H^n}{\Delta t} - M_C \frac{h^{n+1} - h^n}{\Delta t}.
\end{align*}
Here, $\Delta t$ denotes the time increment. We employ the optimization strategy described in \cite{ruppenthal2024scalable}. Introducing the lumped mass matrix allows us to rewrite the discretized equation as
\begin{align} \label{eq:bath_evolution}
     M_L \frac{b^{n+1} - b^n}{\Delta t} = M_C \frac{H^{n+1} - H^n}{\Delta t} - M_C \frac{h^{n+1} - h^n}{\Delta t} + \left( M_L - M_C \right) p.
\end{align}
This formulation is equivalent to the previous one for the particular choice
\begin{align*}
    p = \frac{b^{n+1} - b^n}{\Delta t}.
\end{align*}
In general, the adjustable term $(M_L-M_C) p$ should regularize the inverse problem and prevent spurious bathymetry modes due to possible non-uniqueness. As noticed in \citep{hajduk2020bathymetry}, it resembles the Brezzi-Pitkäranta pressure stabilization term for the incompressible Navier-Stokes equations. The graph Laplacian operator $M_L-M_C$ was used for this purpose by Becker and Hansbo \cite{becker2008simple}. This matrix has zero row sums and a non-trivial kernel spanned by the vector $\mathbbm 1$. The vector $p$ is composed from \emph{flux potentials} that can serve as optimization variables, as shown in \cite{ruppenthal2024scalable}.

We introduce an objective functional that quantifies the discrepancy between the simulated and measured free surface elevation
\begin{align} \label{eq:J1}
    J_1(b,p) = \frac{\alpha}{2} \left\|h^{n+1} + b - H^{n+1}\right\|^2_{M_L} + \frac{\beta}{2} \left\|p\right\|^2_{M_L} + \frac{\gamma}{2} \left\| b - b_e \right\|^2_{\bmat}.
\end{align}
In this context, the norm $\left\|\cdot\right\|_{M_L}$ is induced by the lumped mass matrix
\begin{align*}
    \left\| p \right\|_{M_L} = \sqrt{p^\top M_L p},
\end{align*}
while $\left\| \cdot \right\|_{\bmat}$ denotes the $L^2$ norm on $\Gamma$ induced by the lumped boundary mass matrix $\bmat$. The function $b_e$ incorporates known bathymetry information on $\Gamma$.
The first component of $J_1(b,p)$ penalizes deviations from the measured data $H$, whereas the remaining part introduces a Tikhonov regularization that guarantees uniqueness of optimal flux potentials $p$.
The positive weights $\alpha$, $\beta$, and $\gamma$ determine the relative importance of the corresponding penalty and regularization terms. The use of $\gamma>0$ is essential for ensuring uniqueness of the reconstruction, as it imposes the Dirichlet boundary condition on the bathymetry $b$ in a weak sense and removes the additive null space.

We consider the following problem
\begin{align*}
    \minimize_{b,p \in\mathbb{R}^{N_h}}\;\; J_1(b,p) \quad\text{subject to \eqref{eq:bath_evolution} and}\quad m_i \frac{u^{n+1}_i - u^n_i}{\Delta t} = \sum_{j\in\mathcal N^*_i}2d_{ij}(\overline{u}_{ij}^{b*}-u_i)\quad\forall i.
\end{align*}
In the discretized system of shallow water equations, which serves as a PDE constraint for optimal control, $u$ stands for ($h$, $\hv$).

The Lagrangian of the above optimization problem is given by
\begin{align*}
    \Lambda\left(b, p, q\right) &= \frac{\alpha}{2} \left\|h^{n+1} + b - H^{n+1}\right\|^2_{M_L} + \frac{\beta}{2} \left\|p\right\|^2_{M_L} + \frac{\gamma}{2} \left\| b - b_e \right\|^2_{\bmat} \\[0.25cm]
    &\hphantom{=}\; + \Delta t q^\top \left[ M_L \frac{b - b^n}{\Delta t} - M_C \left( \frac{H^{n+1} - H^n}{\Delta t} - \frac{h^{n+1} - h^n}{\Delta t} \right) - \left( M_L-M_C \right) p \right].
\end{align*}
To circumvent an explicit division by $\Delta t$, we multiply the state equation by $\Delta t$. After this scaling, we proceed to derive the first-order optimality (Karush-Kuhn-Tucker) conditions
\begin{align*}
    \pd{\Lambda}{b} &= \alpha M_L \left( h^{n+1} + b - H^{n+1} \right) + \gamma \bmat \left(b - b_e \right) + M_L q \overset{!}{=} 0, \\
    \pd{\Lambda}{p} &= \beta M_L p - \Delta t \left( M_L - M_C \right) q \overset{!}{=} 0, \\
    \pd{\Lambda}{q} &= M_L \left(b - b^n \right) - M_C \left( H^{n+1} - H^n - h^{n+1} + h^n \right) - \Delta t \left( M_L-M_C \right) p \overset{!}{=} 0.
\end{align*}
This leads to the linear system
\begin{align*}
    \begin{bmatrix} \alpha M_L + \gamma \bmat & 0 & M_L \\ 0 & \beta M_L & -\Delta t (M_L-M_C) \\ M_L & -\Delta t (M_L-M_C) & 0 \end{bmatrix} \begin{bmatrix} b\\p\\q \end{bmatrix} = \mathbf z
\end{align*}
with the right-hand side
\begin{align} \label{eq:ls_rhs}
    \mathbf z =
    \begin{bmatrix}
     \alpha M_L \left( H^{n+1} - h^{n+1} \right) + \gamma \bmat b_e \\
     0 \\
     M_L b^n + M_C \left( H^{n+1} - H^n - h^{n+1} + h^n \right)
    \end{bmatrix}.
\end{align}
Notably, the variable $b$ can be eliminated from this system. After this manipulation, the resulting reduced system takes the form
\begin{align} \label{eq:reduced_system}
    \begin{bmatrix} \beta M_L & -\Delta t (M_L-M_C) \\ -\Delta t (M_L-M_C) & - M_L \left(\alpha M_L + \gamma \bmat \right)^{-1} M_L \end{bmatrix} \begin{bmatrix} p \\ q \end{bmatrix} = \begin{bmatrix} 0 \\ z \end{bmatrix}
\end{align}
with the reduced right-hand side
\begin{align*}
    z &= M_L b^n + M_C \left( H^{n+1} - H^n - h^{n+1} + h^n \right) \\
    &\hphantom{=}\;+ M_L \left(\alpha M_L + \gamma \bmat \right)^{-1} \left( \alpha M_L \left(h^{n+1} - H^{n+1}\right) - \gamma \bmat b_e \right).
\end{align*}
Uniqueness of the flux potentials $p$ and of the Lagrange multiplier $q$ is guaranteed for $\alpha, \beta > 0$. The use of $\gamma > 0$ in the regularization term is required to recover a unique bathymetry~$b$.

\subsection{Reduced-space optimization}\label{sec:sol}\noindent
The optimization process can be made more efficient by incorporating a solution operator $S$ defined so that $S(p) = b$ is obtained through the solution of \eqref{eq:bath_evolution}, which is equivalent to
\begin{gather*}
    M_L \frac{b - b^n}{\Delta t} = M_C \left(\frac{H^{n+1} - H^n}{\Delta t} - \frac{h^{n+1} - h^n}{\Delta t}\right) + \left( M_L - M_C \right) p \\[1em]
    \Leftrightarrow b = M^{-1}_L \left( M_L b^n + M_C \left(H^{n+1} - H^n - h^{n+1} + h^n \right) + \Delta t \left( M_L - M_C \right) p \right).
\end{gather*}
We can express $b$ as an affine linear mapping of the control variable $p$
\begin{align} \label{eq:solution_operator}
    b = T p + r,
\end{align}
where 
\begin{align} \label{eq:T_and_r}
    T = \Delta t M_L^{-1} \left( M_L - M_C \right) \;\text{and}\; r = b^n + M^{-1}_L M_C \left(H^{n+1} - H^n - h^{n+1} + h^n \right).
\end{align}
By eliminating the variable $b$ from the objective given by \eqref{eq:J1}, we obtain an unconstrained optimization problem of the form
\begin{align*}
    \minimize_{p \in \mathbb R^{N_h}}\;\; \tilde J_1(p) = \frac{\alpha}{2} \left\|h^{n+1} + S(p) - H^{n+1}\right\|^2_{M_L} + \frac{\beta}{2} \left\|p\right\|^2_{M_L} + \frac{\gamma}{2} \left\| S(p) - b_e \right\|^2_{\bmat}.
\end{align*}
This reduced-space formulation is equivalent to the reduced system \eqref{eq:reduced_system}.

\section{Regularization strategies}\label{sec:regularization}\noindent
A further challenge is the presence of noise in the measurement data. We augment the optimal-control functional with sparsity-promoting terms that damp non-physical oscillations while still allowing sharp bathymetric features. In this context, sparsity implies that the recovered bathymetry is locally constant with isolated discontinuities. Equivalently, its gradient vanishes almost everywhere, being non-zero only on a set of small measure. Total variation denoising (TVD) regularizes the objective function by penalizing the Euclidean norm of the gradient, which suppresses high-frequency noise while capturing discontinuities. In this section, we consider and discuss multiple $L^1$ gradient regularizations.

\subsection{Total variation denoising \tvd}\noindent
To enhance the reconstruction fidelity, we incorporate the smooth TVD regularizer that was used to remove noise in \cite{vogel1996iterative}. With this additional term, the objective function becomes
\begin{align*}
    J_2(b,p) = \frac{\alpha}{2} \left\|h^{n+1} + b - H^{n+1}\right\|^2_{M_L} + \frac{\beta}{2} \left\|p\right\|^2_{M_L} + \frac{\gamma}{2} \left\| b - b_e \right\|^2_{\bmat} + \epsilon \int_{\Omega} \sqrt{\left | \nabla b \right |^2_2 + \zeta^2 } \dx.
\end{align*}
The regularization is governed by the parameters $\epsilon$ and $\zeta$. If $\epsilon$ is too small, the noise persists, whereas an overly large $\epsilon$ leads to a diffusion-dominated scheme that produces a smeared solution.

The new definition of the objective function yields a nonlinear system of the form
\begin{align*}
    \begin{bmatrix} \alpha M_L + \gamma \bmat + W(b) & 0 & M_L \\ 0 & \beta M_L & -\Delta t (M_L-M_C) \\ M_L & -\Delta t (M_L-M_C) & 0 \end{bmatrix} \begin{bmatrix} b\\p\\q \end{bmatrix} = \mathbf z
\end{align*}
with the right-hand side $\mathbf z$ defined as in \eqref{eq:ls_rhs} and $W(b) = \left( w_{ij} (b) \right)_{i,j=1}^{N_h}$ with
\begin{align*}
    w_{ij}(b) = \epsilon \int_{\Omega} \frac{1}{\sqrt{ \left | \nabla b \right |^2_2 + \zeta^2 }} \nabla \varphi_i \cdot \nabla \varphi_j \dx.
\end{align*}
At first glance, this system cannot be simplified as before and remains nonlinear, which considerably increases the computational effort. Moreover, we now need to recover the bathymetry gradient and assemble a diffusion matrix $W$. The dependence of $W$ on $b$ can be handled either explicitly or implicitly. In the latter case, the matrix and corresponding preconditioners must be re-assembled at every Newton iteration. In our experience, the implicit strategy yields more accurate solutions, whereas the explicit approach using a factorization of $W$ is faster but more diffusive.

\subsection{\Lonelone~ gradient regularization}\noindent
Our numerical experiments have shown that reducing noise in the reconstructed bathymetry is most effectively achieved by incorporating gradient-based regularization. To extend this idea, we augment the objective functional with an $L^1$-regularization term on $\nabla b$. The aim is to achieve sparsity in the bathymetry gradient in order to suppress spurious oscillations arising from measurement noise. As noted by Boyd and Vandenberghe \citep{boyd2004convex}, $L^1$ regularization is known for its feature-selection properties, yielding sparse solutions such that small gradient components are driven towards zero in our case. However, formulating the $L^1$ norm as an integral makes the objective non-differentiable, thereby complicating the numerical solution of the resulting optimization problem. To address this challenge, we employ a dual formulation of the $L^1$ norm. The dual formulation yields a smooth objective function, subject solely to box constraints, and can be solved efficiently by standard optimization algorithms. This formulation is based on the Fenchel-Rockafellar duality argument (cf. \cite{rockafellar1970convex,bauschke2017convex,kunisch2004total}).

Our objective is to minimize the $L^1$ norm
\begin{align*}
    \|\nabla b\|_{L^1(\Omega;\mathbb R^d)} = \int_{\Omega} \left|\nabla b \!\left(\mathbf x \right)\right|_1 \dx,
\end{align*}
where $|\cdot |_1$ denotes the $l_1$ norm on $\mathbb R^d$. The discrete bathymetry gradient $\nabla b_h \left( \mathbf x \right) =\sum_{i=1}^{N_h} b_i\,\nabla\varphi_i \left(\mathbf x \right)$ is discontinuous across element interfaces and piecewise-linear (for quadrilaterals, or piecewise-constant for triangles). Using the dual representation of the $L^1$ norm, we find that
\begin{align}
  \label{eq:dual_L1L1}
  \|\nabla b_h\|_{L^1(\Omega;\mathbb R^d)} = \sup_{\substack{\\[.3em]\mathbf g\in L^{\infty}(\Omega;\mathbb R^d)\\[.3em]|\mathbf g (\mathbf x)|_{\infty}\le 1 \text{ a.e. } \mathbf x \in \Omega}}
    \int_{\Omega}\nabla b_h\! \left(\mathbf x \right)\cdot\mathbf g \!\left(\mathbf x \right) \dx,
\end{align}
where the $l_\infty$ norm represents the dual $l_1$ norm. Applying a quadrature rule
\begin{align} \label{eq:quad_rule}
\mathcal Q = \left\{ \left(\mathbf z_j,\omega_j \right)\mid j=1,\dots,M_h \right\},
\qquad
\mathbf z_j \in\mathbb R^d,\;\;\omega_j > 0
\end{align}
to \eqref{eq:dual_L1L1} yields the discrete optimization problem
\begin{subequations}\label{eq:discrete_opt_l1l1}
\begin{align}
  \maximize_{\{\mathbf g_j\}_{j=1}^{M_h}\in\mathbb R^{dM_h}} &\;\;\sum_{j=1}^{M_h}\omega_j\,\nabla b_h \!\left( \mathbf z_j \right)\cdot\mathbf g_j \\[4pt]
  \text{subject to}\quad&|\mathbf g_j|_\infty \le 1, \qquad j=1,\dots,M_h . \label{eq:discrete_opt_l1_cond}
\end{align}
\end{subequations}
In this case, condition \eqref{eq:discrete_opt_l1_cond} can be replaced by the box constraints
\begin{align*}
    -1 \le \mathbf g_j \le 1, \qquad j=1,\dots,M_h.
\end{align*}
These constraints are again meant to hold component-wise. Extrema of bilinear finite element approximations are attained at the vertices of mesh cells. Consequently, to suppress oscillations it is advantageous to employ a quadrature rule that places its points along element edges, thereby fully capturing the dynamic behavior. For instance, one may use a quadrature rule whose integration points include the nodes of a bilinear discontinuous Galerkin (DG) space on quadrilaterals. Since bilinear Lagrange basis functions possess the interpolation property, we may identify the discrete dual variables $\mathbf g_j$ with values of a vector-valued DG function
\begin{align*}
    \mathbf g_h \!\left(\mathbf z_j\right)=\sum_{i=1}^{M_h}\mathbf g_i\,\psi_i\!\left(\mathbf z_j\right) = \mathbf g_j .
\end{align*}
Thus we can interpret the dual variable $\mathbf g_h$ as a discontinuous, piecewise-linear finite element function and use its values at the quadrature points directly in \eqref{eq:discrete_opt_l1l1}.
To begin, we express the objective functional in the following form
\begin{align} \label{eq:obj3}
     J_3(b,p) \;=\; \frac{\alpha}{2} \left\|h^{n+1} + b - H^{n+1}\right\|^2_{M_L} + \frac{\beta}{2} \left\|p\right\|^2_{M_L} + \frac{\gamma}{2} \left\| b - b_e \right\|^2_{\bmat} + \kappa \left\|\nabla b \right\|_{L^1(\Omega;\mathbb R^d)}.
\end{align}
Next, we eliminate the state variable $b$ by invoking the solution operator $S$ introduced in \Cref{sec:sol}. This yields
\begin{align*}
     \tilde J_3(p) \;=\; \frac{\alpha}{2} \left\|h^{n+1} + S(p) - H^{n+1}\right\|^2_{M_L} + \frac{\beta}{2} \left\|p\right\|^2_{M_L} + \frac{\gamma}{2} \left\| S(p) - b_e \right\|^2_{\bmat} + \kappa \left\| \nabla S(p) \right\|_{L^1(\Omega;\mathbb R^d)}.
\end{align*}
Replacing the $L^1$ norm by the objective function of the discrete optimization problem \eqref{eq:discrete_opt_l1l1} transforms the original minimization problem into the equivalent primal-dual form
\begin{gather*}
    \min_{p \in \mathbb R^{N_h}} \max_{\mathbf{g} \in \mathbb R^{dM_h}}\;\; \frac{\alpha}{2} \left\| h^{n+1} + S(p) - H^{n+1} \right\|^2_{M_L} + \frac{\beta}{2} \left\|p \right\|^2_{M_L} + \frac{\gamma}{2} \left\| S(p) - b_e \right\|^2_{\bmat} + \mathbf{g}^{\top} \mathbf{A}S(p)\\[1em]
    \text{subject to}\quad -\kappa \leq \mathbf g \leq \kappa.
\end{gather*}
Since $\nabla b_h$ is a discontinuous piecewise-linear function, we assemble the gradient matrix
\begin{align*}
    \mathbf A = \left( \mathbf a_{ij} \right)_{\substack{i=1,\dots,M_h\\ j=1,\dots,N_h}}
\end{align*}
as a mixed finite element operator by utilizing continuous linear basis functions $\left( \varphi_j \right)^{N_h}_{j=1}$ in the trial space and discontinuous linear basis functions $\left( \psi_i \right)^{M_h}_{i=1}$ in the test space. The entries of $\mathbf A$ are given by
\begin{align*}
    \mathbf a_{ij}\;=\;\int_{\Omega} \psi_i\,\nabla \varphi_j \dx.
\end{align*}
The dual formulation is well-defined: it is a linear program subject only to box constraints, and thus admits optimal solution(s). The primal part is strictly convex in $p$, guaranteeing uniqueness of the minimizer. By Sion's minimax theorem \cite{sion1958general} we may interchange the order of the maximization and minimization without altering the solution of the problem
\begin{gather*}
    \max_{\mathbf{g} \in \mathbb R^{dM_h}} \min_{p \in \mathbb R^{N_h}} \;\; J_4 (p, \mathbf g) = \frac{\alpha}{2} \left\| h^{n+1} + S(p) - H^{n+1} \right\|^2_{M_L} + \frac{\beta}{2} \left\|p \right\|^2_{M_L} + \frac{\gamma}{2} \left\| S(p) - b_e \right\|^2_{\bmat} + \mathbf{g}^{\top} \mathbf{A}S(p)\\[1em]
    \text{subject to}\quad -\kappa \leq \mathbf g \leq \kappa.
\end{gather*}
Next, we embed the reduced framework introduced in \eqref{eq:solution_operator} for the state $p$ and the control $\mathbf g$ by utilizing equation \eqref{eq:T_and_r} for $T$ and $r$. Since our focus is on the dual variable $\mathbf g$, we eliminate $p$ by deriving its first-order optimality condition
\begin{align*}
   \pd{J_4 (p, \mathbf g)}{p} = \alpha T^{\top} M_L \left( h^{n+1} + S(p) - H^{n+1} \right) + \beta M_L p + \gamma T^{\top} \bmat \left(S(p) - b_e \right) + \left( \mathbf{A}T \right)^{\top} \mathbf g \overset{!}{=} 0.
\end{align*}
We then regroup the terms
\begin{align*}
   \left( \alpha T^{\top} M_L T + \beta M_L + \gamma T^{\top} \bmat T \right) p = \alpha T^{\top} M_L \left( H^{n+1} - h^{n+1} - r \right) - \left( \mathbf{A}T \right)^{\top} \mathbf g + \gamma T^{\top} \bmat \left(b_e - r \right)
\end{align*}
and express $p$ as a function of $\mathbf g$
\begin{align*}
   p( \mathbf g) = \left( T^{\top} \left( \alpha M_L + \gamma \bmat \right) T + \beta M_L \right)^{-1} \left[ T^{\top} \left( \alpha M_L \left( H^{n+1} - h^{n+1} - r \right) - \mathbf{A}^{\top} \mathbf g + \gamma \bmat \left( b_e - r \right) \right) \right].
\end{align*}
With this representation for $p$, the objective functional can be rewritten in the following form
\begin{align*}
    \maximize_{\mathbf{g} \in \mathbb R^{dM_h}} \;\;\hphantom{-} \tilde J_4 \left( \mathbf g \right) \quad\text{subject to}\quad -\kappa \leq \mathbf g \leq \kappa.
\end{align*}
Rather than solving for the maximum directly, we minimize the negated objective instead
\begin{align*}
    \minimize_{\mathbf{g} \in \mathbb R^{dM_h}} \;\;- \tilde J_4 \left( \mathbf g \right) \quad\text{subject to}\quad -\kappa \leq \mathbf g \leq \kappa.
\end{align*}
The resulting formulation is differentiable, subject solely to scalar box constraints. Sparse optimization often suffers from ill-posedness and non-uniqueness issues: the dual variable $\mathbf g$ may require many iterations of the optimization solver, or convergence may stall. To alleviate this, we augment the objective $\tilde J_4$ with an $L^2$-regularization term that depends solely on $\mathbf g$. This modification preserves the optimality conditions for $p$, and the objective becomes
\begin{align*}
    \tilde J_5(\mathbf g) = - \tilde J_4(\mathbf g) + \frac{\nu}{2} \left\| \mathbf g \right \|^2_2.
\end{align*}
As $\tilde J_4$ is concave with respect to $\mathbf g$, the negated objective $-\tilde J_4$ is convex. Augmenting this with a squared $L^2$ norm, which is strictly convex, ensures the total objective remains strictly convex, guaranteeing a unique minimizer for $\tilde J_5$.

\subsection{\Loneltwo~ gradient regularization}\noindent
In the preceding section, we adopted the $l_1$ norm as the vector norm for the bathymetry gradient. By optimizing each component separately, we gain finer control over the anisotropic behavior. Consequently, we now turn our attention to an isotropic $L^1$ formulation that employs the $l_2$ norm on the underlying vector space. The Euclidean norm is the natural choice in a multi-dimensional context, and leads to nonlinear constraints for our problem. The dual approach is preferred since it provides smooth derivatives. To develop a unifying framework, we therefore revert to the original $L^1$ norm formulation before incorporating the Euclidean metric. We redefine
\begin{align*}
    \left\|\nabla b\right\|_{L^1(\Omega;\mathbb R^d)} = \int_{\Omega} \left|\nabla b \!\left(\mathbf x \right) \right|_2 \dx,
\end{align*}
where $|\cdot|_2$ denotes the Euclidean $l_2$ norm on $\mathbb R^d$. The modified dual representation of the $L^1$ norm yields
\begin{align} \label{eq:dual_L1L2}
  \|\nabla b_h\|_{L^1(\Omega;\mathbb R^d)} =\sup_{\substack{\\[.3em]\mathbf g\in L^{\infty}(\Omega;\mathbb R^d)\\[.3em]|\mathbf g (\mathbf x)|_{2}\le 1 \text{ a.e. } \mathbf x \in\Omega}}
    \int_{\Omega}\nabla b_h \!\left(\mathbf x \right) \cdot\mathbf g \!\left(\mathbf x \right) \dx.
\end{align}
Applying the quadrature rule \eqref{eq:quad_rule} to equation \eqref{eq:dual_L1L2} results in the discrete optimization problem
\begin{subequations}
\begin{align*}
  \maximize_{\{\mathbf g_j\}_{j=1}^{M_h}\in\mathbb R^{dM_h}} &\;\;\sum_{j=1}^{M_h}\omega_j\,\nabla b_h \!\left(\mathbf z_j\right)\cdot\mathbf g_j \\[4pt]
  \text{subject to}\quad &|\mathbf g_j|_2 \le 1, \qquad j=1,\dots,M_h .
\end{align*}
\end{subequations}
Recall that we defined the objective function $J_3$, given by \eqref{eq:obj3}. Replacing the $l_1$ norm in the objective with the $l_2$ formulation, we arrive at an equivalent primal-dual problem
\begin{align} \label{eq:prox_method}
    \min_{p \in \mathbb R^{N_h}} \max_{\mathbf{g} \in \mathbb R^{dM_h}}\;\; \frac{\alpha}{2} \left\|S(p) + h^{n+1} - H^{n+1} \right\|^2_{M_L} + \frac{\beta}{2} \left\|p \right\|^2_{M_L} + \frac{\gamma}{2} \left\| S(p) - b_e \right\|^2_{M_\Gamma} + \mathbf{g}^{\top} \mathbf{A}S(p) + \iota_{\mathcal C}(\mathbf{g}),
\end{align}
where the indicator function
\begin{align*}
    \iota_{\mathcal C} (\mathbf{g}) = \begin{cases} 0,& \mathbf{g}\in\mathcal{C},\\ +\infty,&\text{otherwise} \end{cases}
\end{align*}
corresponds to the constraint set
\begin{align*}
    \mathcal{C}= \left\{ \mathbf{g} \in \mathbb{R}^{dM_h} \,\middle|\, \left| \mathbf{g}_j \right|_2 \le \kappa,\, j=1,\dots,M_h \right\}.
\end{align*}
For proximal algorithms, the proximity operator of $\iota_{\mathcal C}$ reduces to a projection \citep{bauschke2017convex} onto the Euclidean ball of radius $\kappa$
\begin{align*}
    \operatorname{Prox}_{\iota_{\mathcal C}}(\mathbf g_j) = \begin{cases}
        \mathbf g_j,&\left| \mathbf g_j \right|_2\leq\kappa, \\[4pt]
        \dfrac{\kappa}{\left| \mathbf g_j \right|_2} \mathbf g_j, &\text{otherwise}.
    \end{cases}
\end{align*}
This projection is applied to all vectors $\mathbf g_j$.

\section{Numerical solution strategy and implementation}\label{sec:implementation}\noindent
Most approaches considered in this work are formulated as minimization problems subject to linear constraints, such as
\begin{gather*}
    \minimize_{b, p \in \mathbb{R}^{N_h}}\;\; J_1 \left( b,p \right) \quad\text{subject to}\quad S(p) = b.
\end{gather*}
Explicitly optimizing the state variable $b$ increases the system dimensionality relative to the vector of flux potentials $p$. To mitigate this growth in problem size, we employ the SimOpt interface provided by ROL. The SimOpt interface eliminates the state variable $b$ by expressing it as a function of $p$, leading to an unconstrained optimization problem involving only the solution operator $S$
\begin{align*}
    \minimize_{p \in \mathbb{R}^{N_h}}\;\; \tilde J_1 \left( p \right).
\end{align*}
By \eqref{eq:solution_operator}, $S$ is an affine linear mapping and can be evaluated at low cost. However, applying its inverse requires greater computational effort. Because the underlying matrix $M_L - M_C$ is a graph Laplacian, we may employ the same solvers as in our previous work \cite{ruppenthal2024scalable}. With the elimination of $b$, the optimization problem depends solely on $p$.

A second class of problems consists of convex optimization subject to box constraints
\begin{gather*}
     \minimize_{\mathbf g \in \mathbb{R}^{dM_h}}\;\; \tilde J_5 \left( \mathbf g \right) \quad\text{subject to}\quad - \kappa \leq \mathbf g \leq \kappa.
\end{gather*}
We apply the same SimOpt interface for the reduced-space formulation. SimOpt constructs a reduced objective incorporating gradient and Hessian information via the chain rule, with the solution operator $S$ integrated into the constraint module. This module performs forward and backward solves between state and control variables. Users supply only the derivatives of the full-size problem.

To solve the bound-constrained or reduced-space problems, we choose a trust-region algorithm designated as `Lin-Moré' in ROL. Building on the trust-region framework introduced in \cite{lin1999newtons}, it seeks to minimize the quadratic model
\begin{align*}
    m_{\mathbf k} \left( \mathbf g \right) = \left( \nabla^2 \tilde J_5 \left( \mathbf g_{\mathbf k} \right) \left( \mathbf g - \mathbf g_{\mathbf k} \right), \mathbf g - \mathbf g_{\mathbf k} \right) + \left(\nabla \tilde J_5 \left( \mathbf g_{\mathbf k} \right), \mathbf g - \mathbf g_{\mathbf k}\right) + \tilde J_5 \left( \mathbf g_{\mathbf k} \right).
\end{align*}
This expression is equivalent to a second-order Taylor series expansion of the objective $\tilde J_5$. The feasible set is defined as
\begin{align*}
    \mathcal C = \left\{ \mathbf g \in \mathbb{R}^{dM_h} \,\middle|\, -\kappa \leq \mathbf g \leq \kappa \right\}.
\end{align*}
During the trust-region iteration, the subproblem solver minimizes a quadratic model subject to both the radius constraint $\tau_{\mathbf k}$ and the feasibility set $\mathcal C$. In the $\mathbf k$-th step this minimization problem reads
\begin{gather*}
    \minimize_{\mathbf g \in \mathbb R^{dM_h}}\;\; m_{\mathbf k} \left( \mathbf g \right) \quad\text{subject to}\quad \left\| \mathbf g - \mathbf{g_k} \right\| \leq \tau_{\mathbf k}, \quad \mathbf g \in \mathcal C.
\end{gather*}
Projection onto the feasible set is performed via component-wise maximum and minimum operations. To solve the optimization problem \eqref{eq:prox_method}, we employ a proximal trust-region framework based on the algorithm proposed by Baraldi and Kouri \citep{baraldi2023proximal}. The resulting subproblems are handled with nonlinear conjugate gradient (NCG) iterations for small-scale instances, while larger problems are addressed using the spectral projected gradient (SPG) method \citep{kouri2022matrix}.

\subsection{Convergence acceleration strategies}\noindent
Several strategies exist to optimize the aforementioned framework. For instance, one could limit changes of bathymetry in individual time steps to improve convergence properties. In the cost functional 
\begin{align*}
    J_6(b,p) = \frac{\alpha}{2} \left\|h^{n+1} + b - H^{n+1}\right\|^2_{M_L} + \frac{\beta}{2} \left\|p\right\|^2_{M_L} + \frac{\gamma}{2} \left\| b - b_e \right\|^2_{\bmat} + \frac{\delta}{2} \left\|b - b^n\right\|^2_{M_L},
\end{align*}
the weight $\delta$ determines how much the bathymetry may vary from one pseudo-time step to the next. If the iterative sequence exhibits oscillations, adjusting $\delta$ allows us to stabilize it similarly to Anderson acceleration \citep{walker2011anderson}. In most theoretical developments, the data is presumed to be sufficiently smooth so that the iterative scheme converges to a steady state. When measurements contain sharp discontinuities or high-frequency components, this assumption is violated: the sequence of iterates $b(\xi)$ may oscillate from one step to the next, and convergence can stall or divergence is possible. A practical remedy is to replace each new iterate by a weighted average with its predecessor, like in approach $J_6$, thereby damping high-frequency noise. More sophisticated schemes, such as Anderson acceleration \cite{walker2011anderson}, construct a linear combination of several past iterates to extrapolate toward the fixed point and have been shown to accelerate convergence in the presence of non-smooth data.

\subsection{Implementation details}\noindent
The optimization problems formulated above are implemented in \verb|C++| using a combination of state-of-the-art scientific libraries. For discretization and linear algebra, we use MFEM (Modular Finite Element Library \cite{anderson2021mfem}), while the optimization routines are supplied by ROL (Rapid Optimization Library \cite{javeed2022get,rol-website}) as part of the Trilinos project \cite{heroux2005overview,heroux2012new,trilinos-website}. MFEM handles mesh generation, finite-element assembly, and data structures. ROL provides a rich set of optimization solvers and interfaces for reduced-space optimization.

All linear systems are solved with either preconditioned GMRES/MINRES iterations or a direct factorization. For distributed factorizations we employ Intel's Cluster Sparse Solver. One-dimensional problems produce small matrices and are therefore handled directly, as this is faster than an iterative approach. In two-dimensional settings of larger scale, we switch to an iterative strategy with multigrid (BoomerAMG) or incomplete LU (Euclid) preconditioners from the \texttt{hypre} library \cite{hypre}. If $\alpha$ and $\gamma$ are comparable, the multigrid version converges rapidly. If $\gamma\gg\alpha$, then the ILU version is preferred because of a possible stall in multigrid convergence. To combine MFEM's efficient backends with ROL's matrix-free framework, we implement a custom vector and operator interface that maps MFEM data structures into ROL objects.

\section{Numerical experiments}\label{sec:experiments}\noindent
We conduct a series of numerical experiments to evaluate the proposed methods. First, we conduct a convergence study on a one-dimensional domain to verify that the algorithms attain the expected order of accuracy and are internally consistent. Next, we extend the problem to two spatial dimensions in order to assess how well the reconstruction procedure handles non-smooth bathymetric features. Finally, we test the robustness of the method by reconstructing the bathymetry from synthetic data contaminated with artificial noise. For all experiments under investigation, we set~$g = \SI{9.81}{\metre\per\second\squared}$. \begingroup \Crefname{section}{Appendix}{Appendices} All parameter values used in the convergence, noise-robustness, and large-scale experiments are listed in \Cref{app:parameter}. \endgroup

\subsection{Convergence study}\label{sec:conv}
\begin{figure}
    \centering
    \begin{subfigure}{0.45\textwidth}
        \centering
        \includegraphics[width=\linewidth,trim={0 4cm 0 4cm},clip]{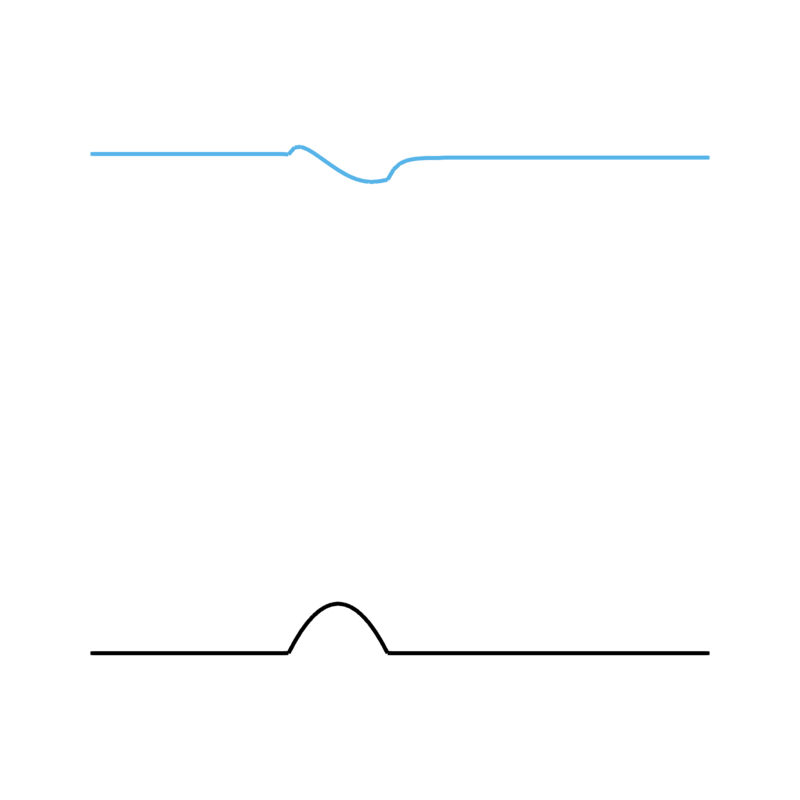}
        \caption{Low-order (ALF) forward solution}
        \label{fig:low_order_solution}
    \end{subfigure}
    \hfill
    \begin{subfigure}{0.45\textwidth}
        \centering
        \includegraphics[width=\linewidth,trim={0 4cm 0 4cm},clip]{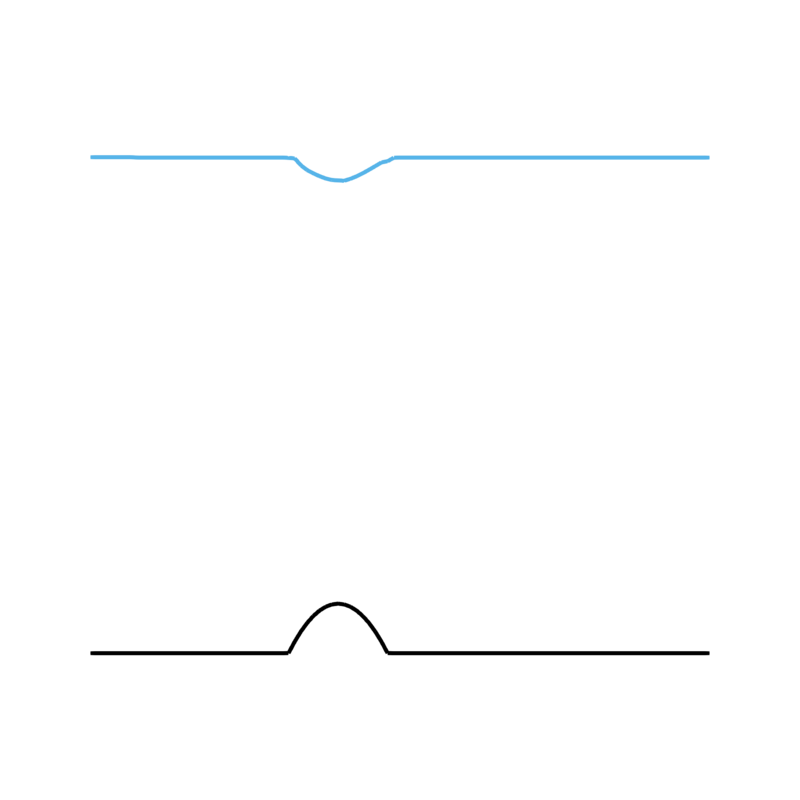}
        \caption{High-order (MCL) forward solution}
        \label{fig:high_order_solution}
    \end{subfigure}
    \caption{Comparison of numerical solutions obtained using the ALF and MCL schemes at the final time $T=\SI{200}{\second}$.}
    \label{fig:solutions}
\end{figure}

\noindent
To verify the consistency of our discretization, we employ a one\babelhyphen{-}dimensional benchmark presented in \cite{vazquezcendon1999improved}. The test domain is a $\SI{25}{\metre}$ long channel discretized with 100 uniformly spaced elements ($\Delta x = \SI{0.25}{\metre}$). The analytical bathymetry is defined as
\begin{align*}
    \tilde b(x) =
    \begin{cases}
        0.2 - 0.05 \,\left(x - 10\right)^2, & \text{if } 8 \leq x \leq 12, \\
        0, & \text{otherwise},
    \end{cases}
\end{align*}
and the initial conditions are $H = \SI{2}{\metre}$ and $v = \SI{2.21}{\metre\per\second}$. The boundary data is chosen consistently with these states. Measurement data is generated by the local Lax-Friedrichs solver, and the simulation is run up to $T=\SI{200}{\second}$ using a time step $\Delta t = \SI{0.03}{\second}$ in conjunction with an SSP-RK2 time integrator.
\Cref{fig:low_order_solution,fig:high_order_solution} display the low-order and limited high-order numerical solutions, respectively. Both figures have been scaled by a factor of $0.1$ in width so that they can be displayed side by side within a single page.

We carry out a convergence test on three successive mesh refinements, with the outcomes listed in \Cref{tab:conv_study_low_order}. We set $\beta = \num{e-11}$ for the optimal control (OC) solution. Since the ALF method provides sufficient diffusion, we also consider a variant with $p=0$ (i.e., no stabilization, NOS) for comparison.

\begin{figure}
    \centering
    \includegraphics[width=0.8\linewidth,trim={0 6cm 0 12cm},clip]{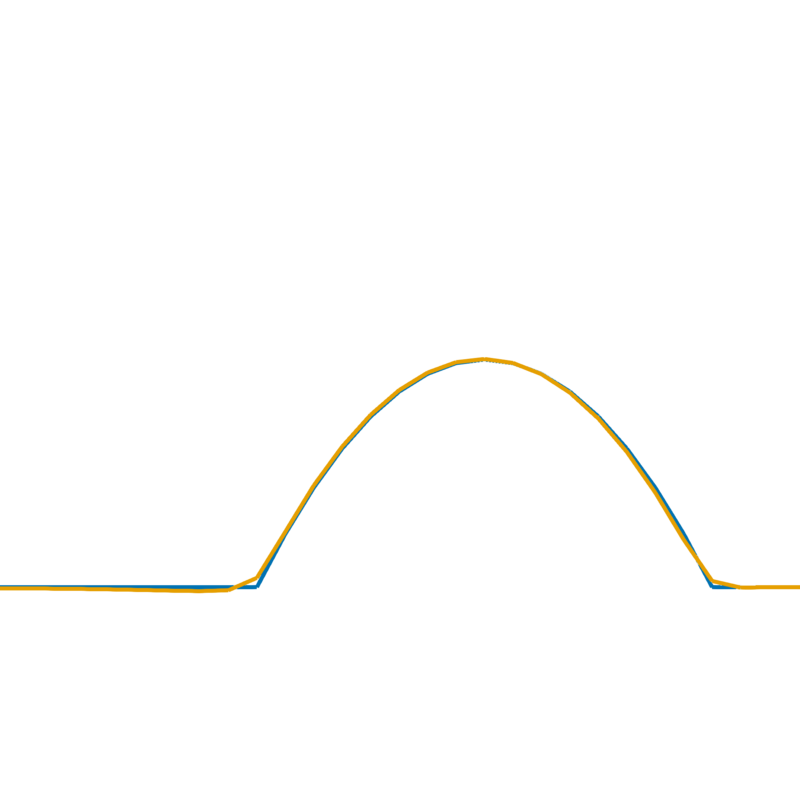}
    \caption{Zoomed view comparing the unstabilized (orange) and optimal control (blue) reconstructions at final time.}
    \label{fig:conv_study_low_order}
\end{figure}

\begin{table}
    \centering
    \begin{tabular}{ccccc}
       Mesh size $\Delta x$ & $L^2$ error (NOS) & Rate (NOS) & $L^2$ error (OC) & Rate (OC) \\
       \midrule
       $1 \slash 4$ & \num{6.50e-3} & --- & \num{1.14e-3} & --- \\[0.1em]
       $1 \slash 8$ & \num{2.45e-3} & $1.41$ & \num{2.85e-4} & $2.00$ \\[0.1em]
       $1 \slash 16$ & \num{9.06e-4} & $1.44$ & \num{7.13e-5} & $2.00$ \\[0.1em]
       $1 \slash 32$ & \num{3.27e-4} & $1.47$ & \num{1.79e-5} & $2.00$ \\
       \bottomrule
    \end{tabular}
    \caption{Convergence rates and $L^2$ errors for bathymetry reconstruction using the ALF scheme (NOS vs. OC).}
    \label{tab:conv_study_low_order}
\end{table}

\Cref{fig:conv_study_low_order} presents the final-time snapshot of two reconstructions: the unstabilized (orange) and the optimal control (blue) reconstruction (objective given by \eqref{eq:J1}). The optimal control reconstruction aligns closely with the reference solution, while the unstabilized approach fails to resolve the sharp corners of the central hump. The convergence rates presented in \Cref{tab:conv_study_low_order} quantify the theoretical expectations for the proposed framework. While the unstabilized scheme exhibits suboptimal convergence ($\approx 1.47$), the optimal control formulation consistently achieves second-order accuracy (rate $2.00$) across all mesh refinements. This confirms that the regularization term effectively mitigates the numerical diffusion inherent in the low-order flux without compromising asymptotic consistency, thereby validating the stability of the discretization strategy for inverse problems.

\begin{figure}
    \centering
    \begin{subfigure}{0.45\textwidth}
        \centering
        \includegraphics[width=\linewidth,trim={0 12cm 0 12cm},clip]{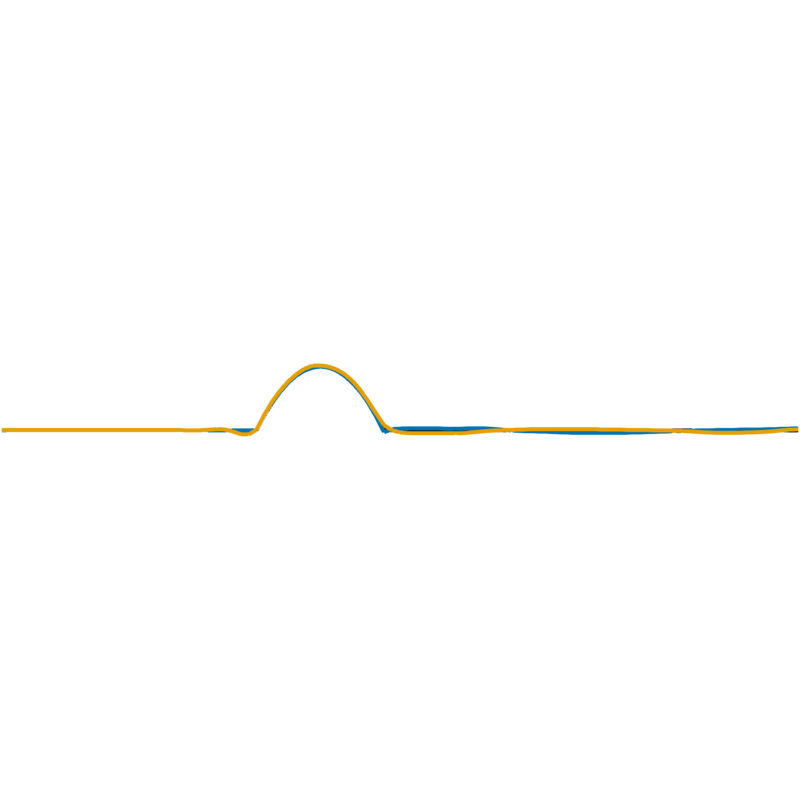}
        \caption{Unstabilized (orange) versus optimal control (blue) comparison}
        \label{fig:high_order_solution_compare}
    \end{subfigure}
    \hfill
    \begin{subfigure}{0.45\textwidth}
        \centering
        \includegraphics[width=\linewidth,trim={0 6cm 0 12cm},clip]{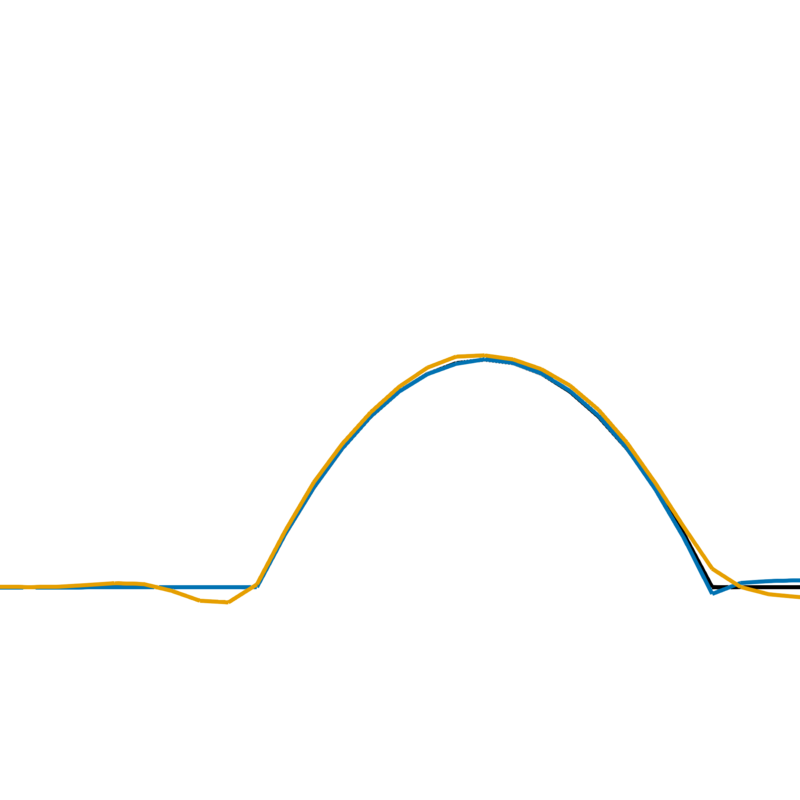}
        \caption{Zoomed view of the hump region including the interpolant (black)}
        \label{fig:high_order_solution_compare_with_zoom}
    \end{subfigure}
    \caption{High-order MCL reconstructions using optimal control regularization ($\beta = \num{e-4}$) at final time $T=\SI{200}{\second}$.}
    \label{fig:high_order_solution_compare_plots}
\end{figure}

We also test a reconstruction with the MCL scheme, using $\beta = \num{e-4}$. Optimal control outperforms the unstabilized approach (see \Cref{fig:high_order_solution_compare_plots}). It differs from the reference at the lower right corner of the hump. In this non-steady configuration, the modified MCL algorithm oscillates between high- and low-order modes, producing small spurious waves, while the unstabilized approach generates even larger oscillations and poorly resolves the left corner of the hump.

\subsection{Two-dimensional benchmark test}\label{sec:two}
\begin{figure}
    \centering
    \includegraphics[width=0.8\linewidth,trim={0 4cm 0 0.4cm},clip]{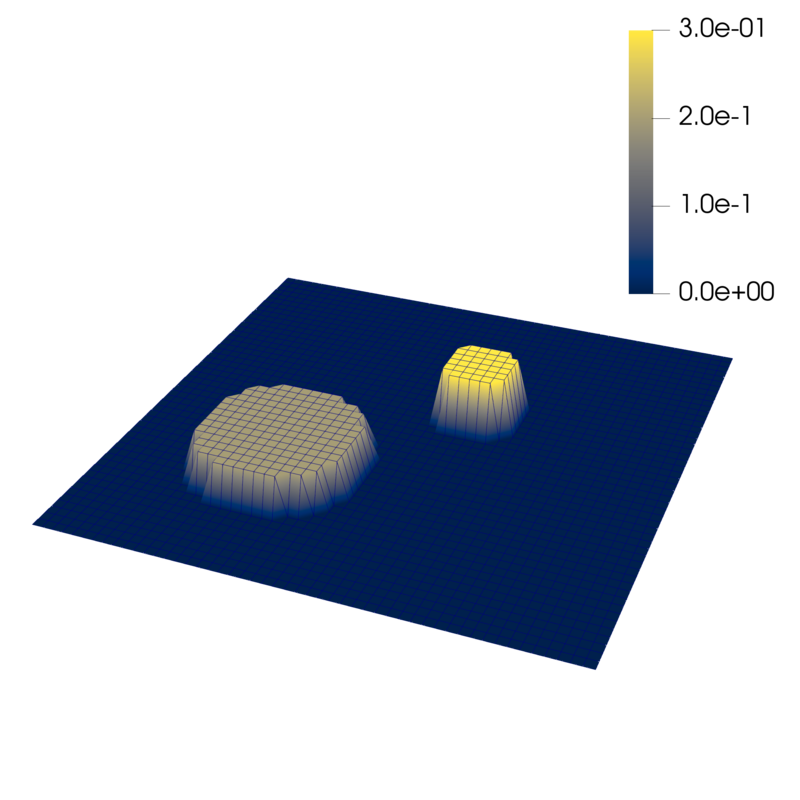}
    \caption{Reference bathymetry profile for the two-dimensional cylinder benchmark.}
    \label{fig:2d_bathymetry}
\end{figure}

\noindent
We extend the one-dimensional benchmark to a square domain of size $\SI{25}{\metre} \times \SI{25}{\metre}$ discretized with a uniform $50 \times 50$ grid. The bathymetry consists of two cylinders
\begin{align*}
    \tilde b(x,y) =
    \begin{cases}
        0.2, & \text{if } \sqrt{\left(x - 8.0 \right)^2+ \left(y - 8.0 \right)^2} \leq 4.0, \\
        0.3, & \text{if } \sqrt{\left(x - 15.0 \right)^2+ \left(y - 15.0 \right)^2} \leq 2.0, \\
        0, & \text{otherwise.}
    \end{cases}
\end{align*}
\Cref{fig:2d_bathymetry} illustrates these cylinders placed sequentially along the flow direction. The first one is lower but wider, while the second is higher and narrower. The initial free surface elevation is $H=\SI{2}{\metre}$ everywhere. The uniform velocity field $\mathbf{v} = \left(2.21, 2.21 \right)^{\top} \si{\metre\per\second}$ is used in this experiment. We apply only the modified MCL scheme to this benchmark. No additional stabilization is used. For optimal control regularization, we set $\beta = \num{e-7}$. The simulation runs until $T=\SI{60}{\second}$ with a time step of $\Delta t = \SI{0.01}{\second}$.

\begin{figure}
    \centering
    \begin{subfigure}{0.45\textwidth}
        \centering
        \includegraphics[width=\linewidth,trim={0 4cm 0 0.4cm},clip]{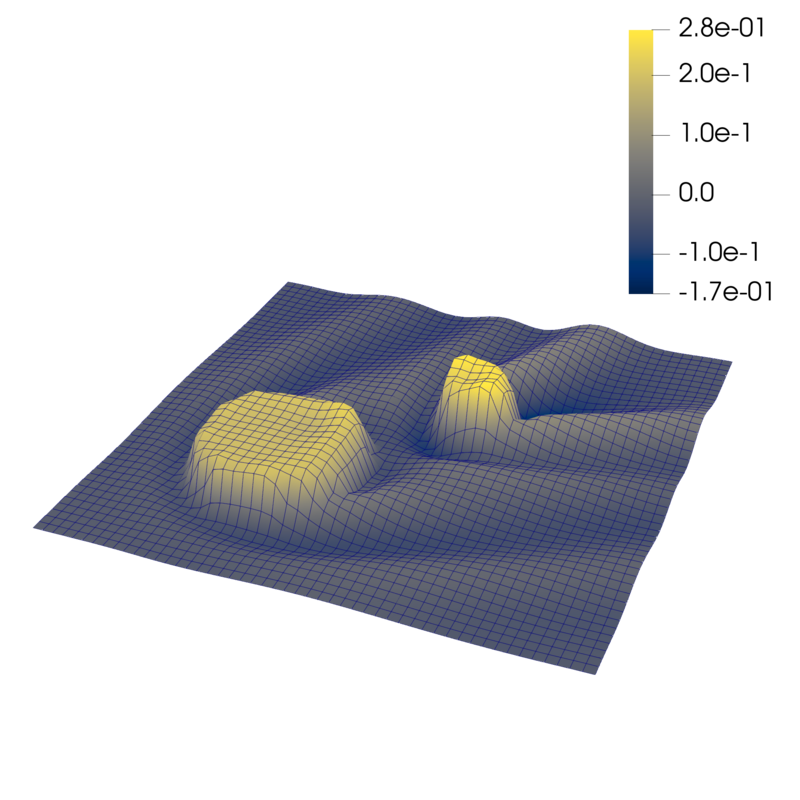}
        \caption{Unstabilized MCL reconstruction ($p=0$)}
        \label{fig:2d_mcl_solution_nos}
    \end{subfigure}
    \hfill
    \begin{subfigure}{0.45\textwidth}
        \centering
        \includegraphics[width=\linewidth,trim={0 4cm 0 0.4cm},clip]{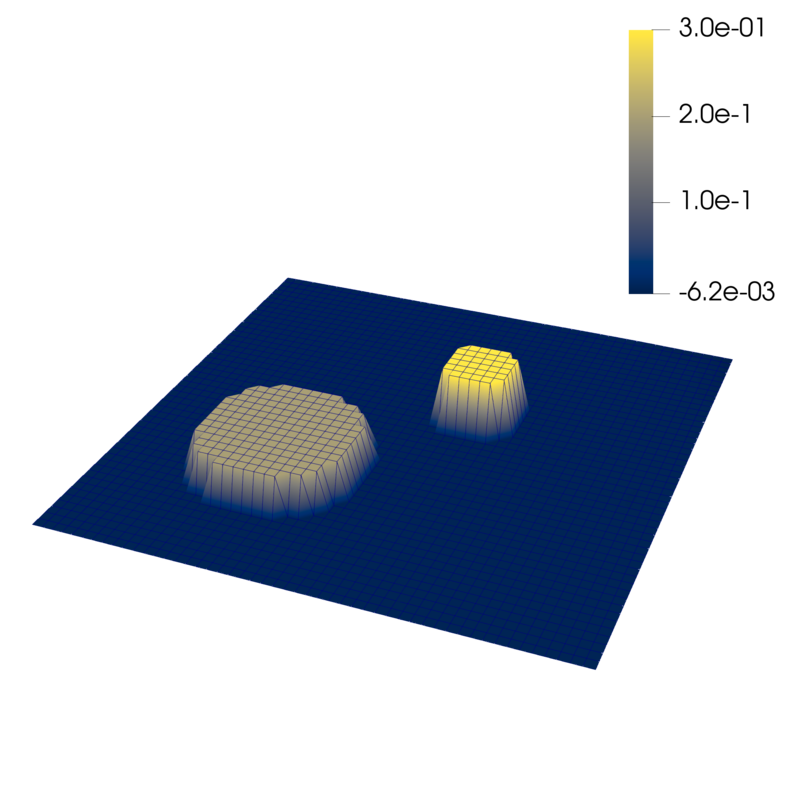}
        \caption{Optimal control regularized MCL reconstruction}
        \label{fig:2d_mcl_solution_oc}
    \end{subfigure}
    \caption{Two-dimensional bathymetry reconstructions after $\SI{60}{\second}$ of simulation time using the MCL scheme.}
    \label{fig:2d_mcl_solutions}
\end{figure}

\Cref{fig:2d_mcl_solutions} shows the two reconstructions. In \Cref{fig:2d_mcl_solution_nos}, the unstabilized reconstruction displays prominent artifacts driven by steep gradients of the cylinders, dominating much of the bathymetry. By contrast, the optimal control reconstruction in \Cref{fig:2d_mcl_solution_oc} lacks these artifacts and closely matches the interpolated bathymetry shown in \Cref{fig:2d_bathymetry}.

\begin{table}
    \centering
    \begin{tabular}{ccc}
       Mesh size $\Delta x$ & $L^2$ error (NOS) & $L^2$ error (OC) \\
       \midrule
       \SI{0.5}{\meter} & \num{9.74e-1} & \num{4.31e-1} \\
       \bottomrule
    \end{tabular}
    \caption{Summary of $L^2$ errors for two-dimensional cylinder benchmark reconstructions using the MCL scheme.}
    \label{tab:two_error}
\end{table}

In the two-dimensional setting, \Cref{tab:two_error} highlights a substantial reduction in reconstruction error when applying optimal control regularization compared to the unstabilized approach. The $L^2$ error decreases by approximately 56\% (from \num{9.74e-1} to \num{4.31e-1}), validating the necessity of stabilization for maintaining solution fidelity in higher dimensions where gradient artifacts are more pronounced. This improvement underscores the efficacy of the flux-potential framework in suppressing spurious oscillations that typically dominate unregularized reconstructions near steep bathymetric features.

\subsection{Noisy data}\label{sec:noisy}\noindent
Up to this point, all test data have been produced by a deterministic forward solver. In practice, however, observations are corrupted by measurement errors, atmospheric disturbances, or other sources of uncertainty. To assess the robustness of our algorithms under such conditions, we generate synthetic noisy data as follows. First, we compute the reference solution and perturb each observation by multiplying with independent Gaussian noise. The perturbation is drawn from $\mathcal N(0,\sigma^2)$, i.e., it has zero mean and a standard deviation $\sigma$ that we vary between $1\%$ and $5\%$ of the signal magnitude. The resulting dataset is used to reconstruct the bathymetry, allowing us to evaluate how well the proposed methods tolerate observational noise. No additional filtering or preprocessing is applied. Any regularization required for the inverse problem is handled within the optimization framework.

\begin{figure}
    \centering
    \begin{subfigure}{0.45\textwidth}
        \centering
        \includegraphics[width=\linewidth,trim={0 9.2cm 0 9.2cm},clip]{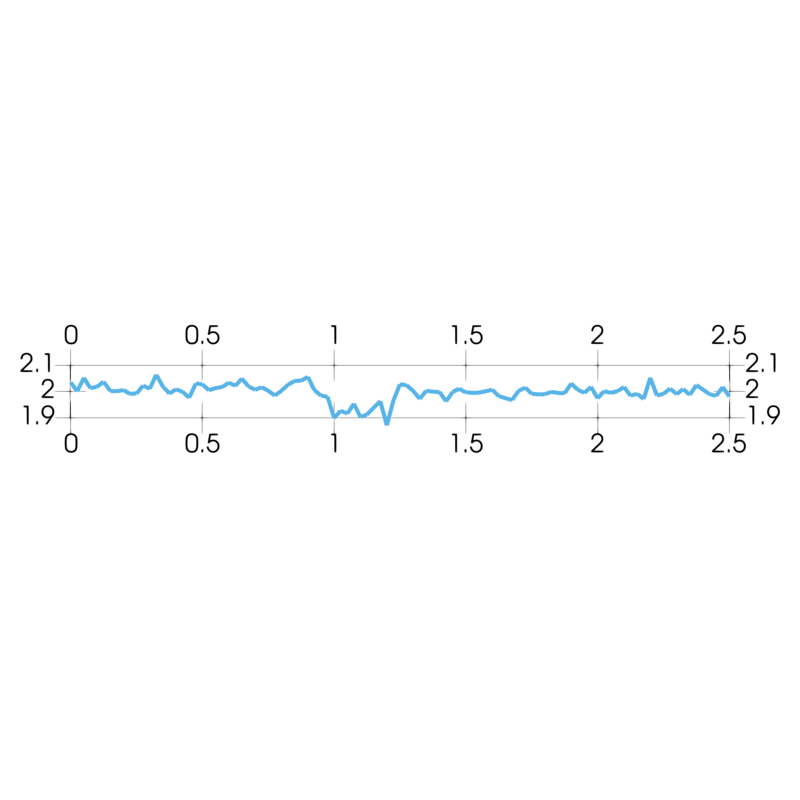}
        \caption{Noisy measurement for $\sigma=1\%$}
        \label{fig:low_order_solution_noise_0.01}
    \end{subfigure}
    \hfill
    \begin{subfigure}{0.45\textwidth}
        \centering
        \includegraphics[width=\linewidth,trim={0 9.2cm 0 9.2cm},clip]{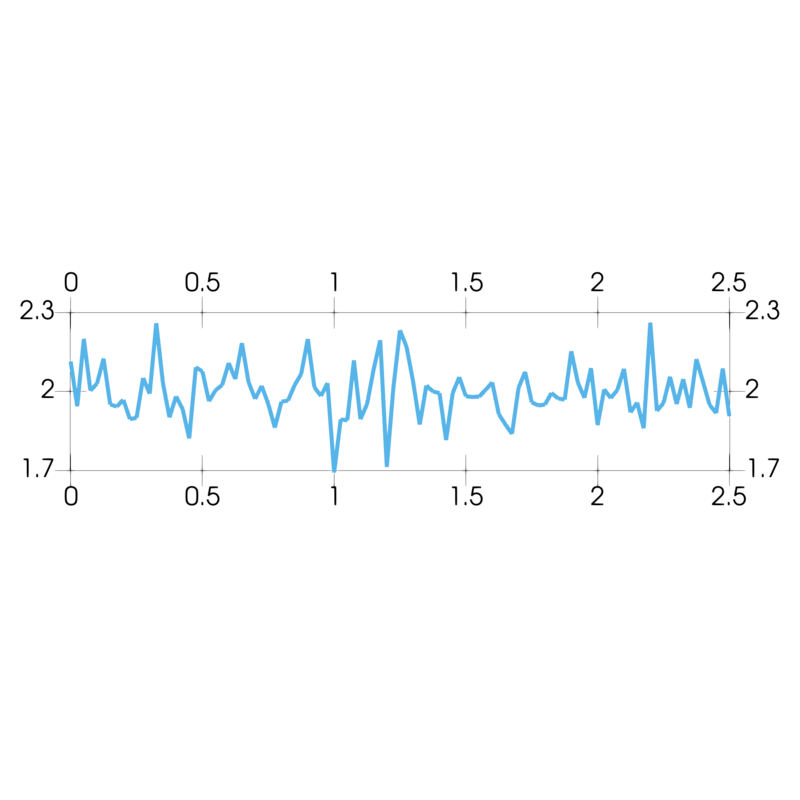}
        \caption{Noisy measurement for $\sigma=5\%$}
        \label{fig:low_order_solution_noise_0.05}
    \end{subfigure}
    \caption{Synthetic noisy observations at final time $T=\SI{200}{\second}$ with standard deviations $\sigma=1\%$ and $\sigma=5\%$.}
    \label{fig:measure_data_with_noise_levels}
\end{figure}

\Cref{fig:measure_data_with_noise_levels} shows two representative noisy data sets for the one-dimensional benchmark case, with standard deviations of $1\%$ and $5\%$. Unlike the clean free-surface simulations in \Cref{fig:solutions}, the characteristic hump becomes obscured by noise.

\begin{figure}
    \centering
    \begin{subfigure}{0.47\textwidth}
        \centering
        \includegraphics[width=\linewidth,trim={0 11.6cm 0 11.6cm},clip]{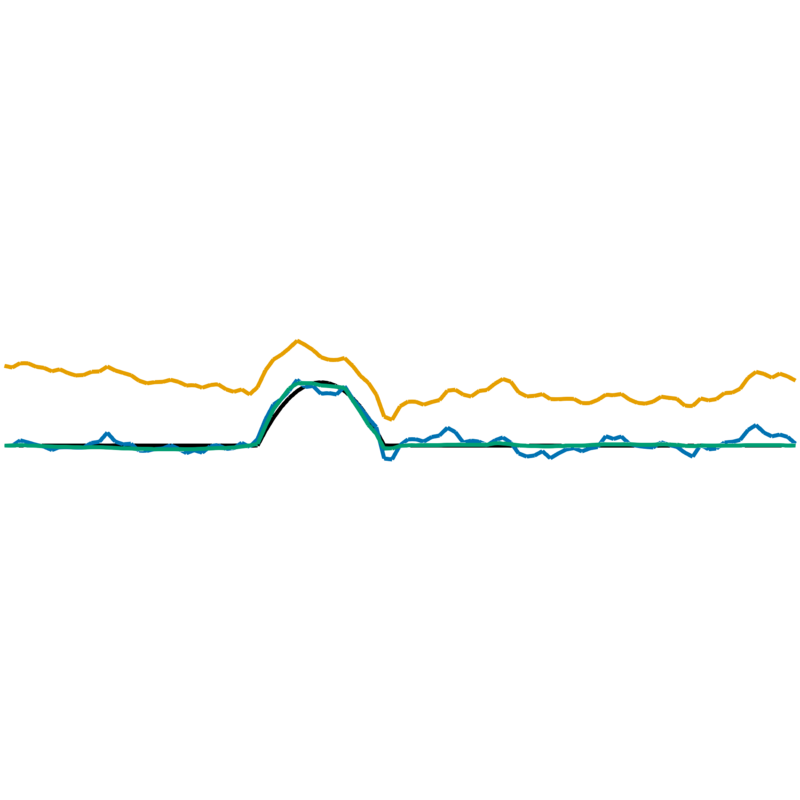}
        \caption{ALF scheme results for $\sigma = 1\%$: unstabilized (NOS, orange), optimal control (OC, blue), and $L^1$ regularization (green)}
        \label{fig:low_order_solution_compare_noise_0.01}
    \end{subfigure}
    \hfill
    \begin{subfigure}{0.47\textwidth}
        \centering
        \includegraphics[width=\linewidth,trim={0 11.6cm 0 11.6cm},clip]{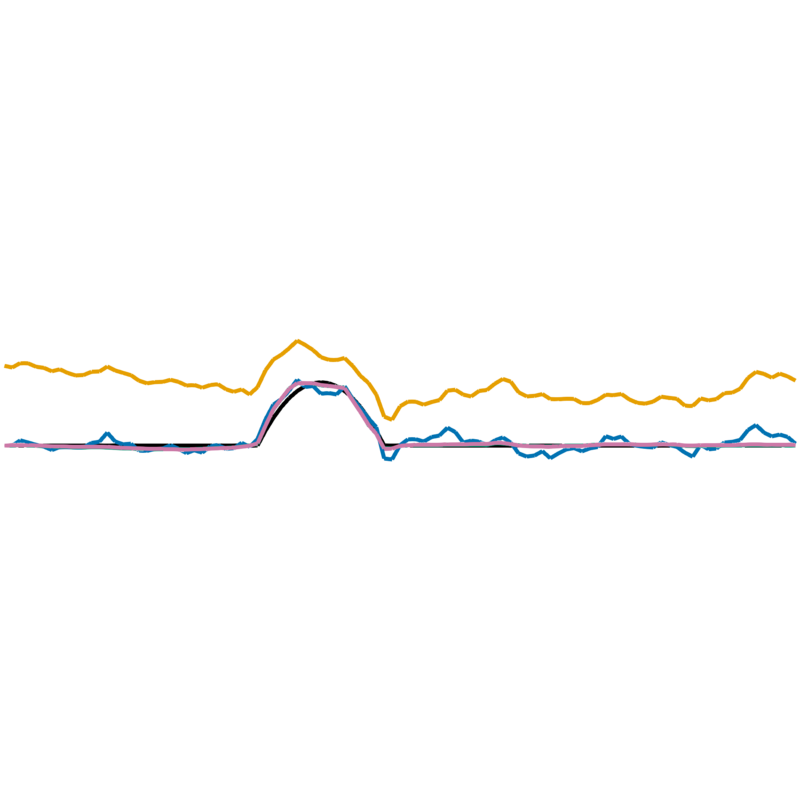}
        \caption{ALF scheme with additional TVD regularization ($\sigma = 1\%$, purple)\\}
        \label{fig:low_order_noise_0.01_l1_l2}
    \end{subfigure}
    \begin{subfigure}{0.47\textwidth}
        \centering
        \includegraphics[width=\linewidth,trim={0 8cm 0 8cm},clip]{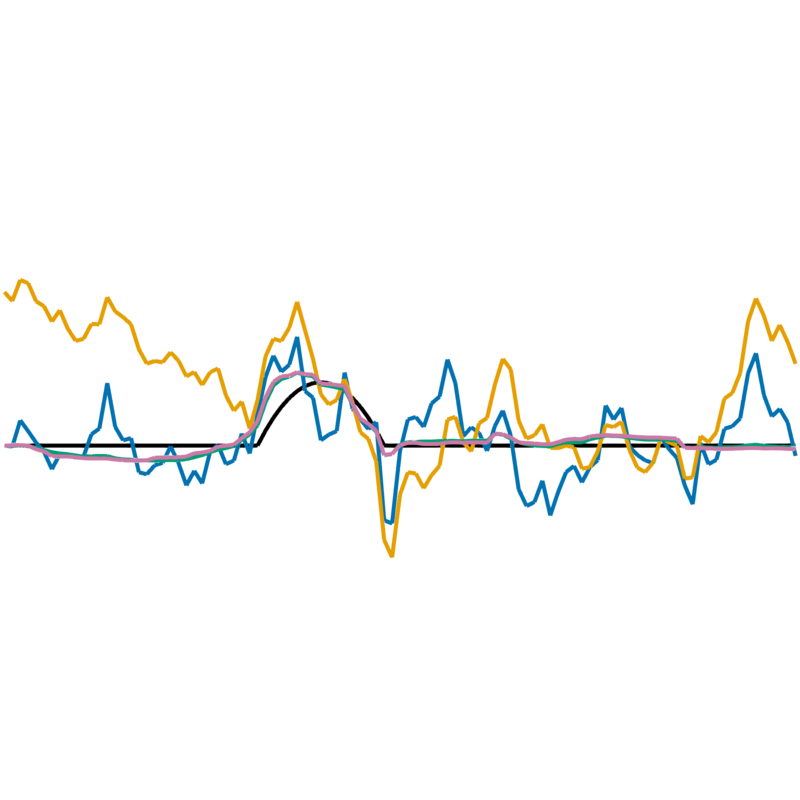}
        \caption{ALF reconstruction for noisy data with $\sigma = 5\%$}
        \label{fig:low_order_solution_compare_noise_0.05}
    \end{subfigure}
    \hfill
    \begin{subfigure}{0.47\textwidth}
        \centering
        \includegraphics[width=\linewidth,trim={0 8cm 0 8cm},clip]{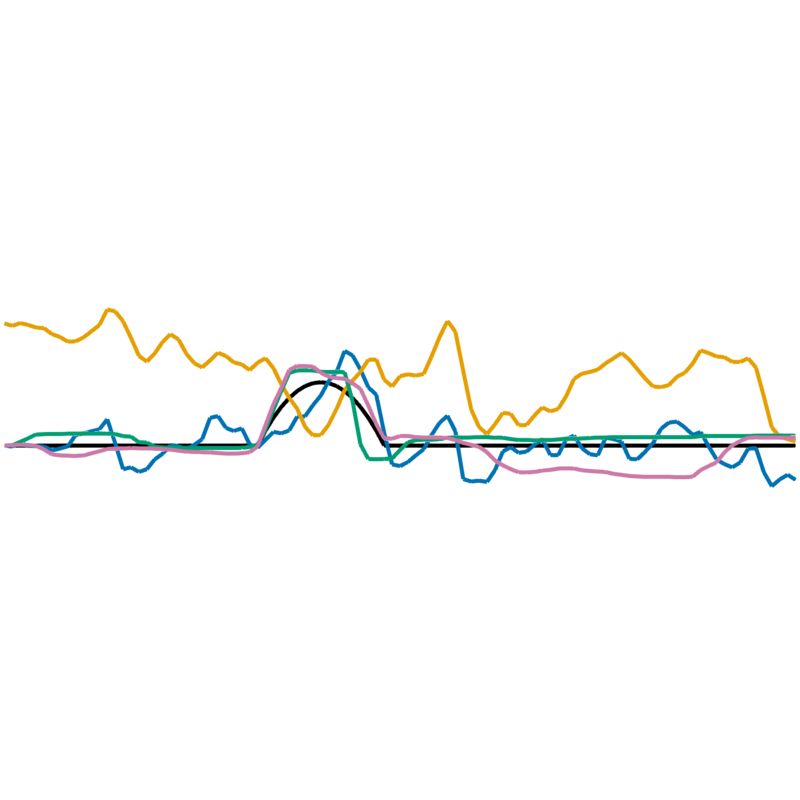}
        \caption{MCL reconstruction for noisy data with $\sigma = 1\%$}
        \label{fig:high_order_noise_0.01}
    \end{subfigure}
    \caption{One-dimensional bathymetry reconstructions under noisy conditions using various regularization strategies.}
    \label{fig:low_order_noise_0.01}
\end{figure}

For a perturbation of $\sigma = 1\%$, both the unstabilized solution (orange, NOS) and the optimal control solution (blue, OC) exhibit pronounced spurious oscillations, see \Cref{fig:low_order_solution_compare_noise_0.01}. The unstabilized approach performs particularly poorly because it cannot enforce boundary conditions, resulting in a large deviation from the true bathymetry.

Adding an $L^1(l_1)$-regularization term (green) flattens the bathymetry reconstruction outside the hump without degrading its shape. \Cref{fig:low_order_noise_0.01_l1_l2} compares this result with a TVD-stabilized reconstruction (purple). The two curves are nearly indistinguishable.

When the noise level is increased to $\sigma=5\%$, the unstabilized and optimal control reconstructions become highly oscillatory, and the hump is no longer recognizable (\Cref{fig:low_order_solution_compare_noise_0.05}). Both the $L^1(l_1)$ regularization and TVD stabilization produce much cleaner estimates that follow the expected bathymetry, although minor disturbances remain near the flat regions and the left side of the hump.

Applying a $1\%$ noise perturbation to the MCL reconstruction (\Cref{fig:high_order_noise_0.01}) yields different behaviors. The TVD stabilization (purple) produces a smooth profile but introduces an artificial dip next to the hump. The $L^1(l_1)$ regularization damps most oscillations but tends to overestimate the steepness on the right side of the hump, resulting in a less smooth profile. When interpreted as an indicator function, the $L^1(l_1)$ solution clearly marks the locations where humps should appear, suggesting it may be preferable for localization tasks.

\begin{table}
    \centering
    \begin{tabular}{ccccc}
       Algorithm ($\sigma$) & $L^2$ error (NOS) & $L^2$ error (OC) & $L^2$ error \tvd & $L^2$ error (\Lonelone) \\
       \midrule
       ALF (1\%) & \num{8.87e-1} & \num{1.05e-1} & \num{4.03e-2} & \num{3.99e-2} \\
       ALF (5\%) & \num{1.19e+0} & \num{5.15e-1} & \num{1.59e-1} & \num{1.48e-1} \\
       MCL (1\%) & \num{1.25e+0} & \num{2.91e-1} & \num{2.57e-1} & \num{1.68e-1} \\
       \bottomrule
    \end{tabular}
    \caption{Summary of $L^2$ errors for one-dimensional noisy data reconstructions under 1\% and 5\% noise levels.}
    \label{tab:noisy_error_1d}
\end{table}

The results in \Cref{tab:noisy_error_1d} demonstrate the robustness of sparsity-promoting regularizers under measurement uncertainty. At $\sigma=1\%$, both TVD and $L^1$ regularization reduce errors by an order of magnitude compared to the unstabilized case, with values dropping from \num{8.87e-1} to approximately \num{4e-2}. Notably, at higher noise levels ($\sigma=5\%$), the $\Lonelone$ variant maintains a slight edge over TVD in preserving feature sharpness while keeping error magnitudes comparable.

\begin{figure}
    \centering
    \includegraphics[width=0.8\linewidth,trim={0 4cm 0 0.4cm},clip]{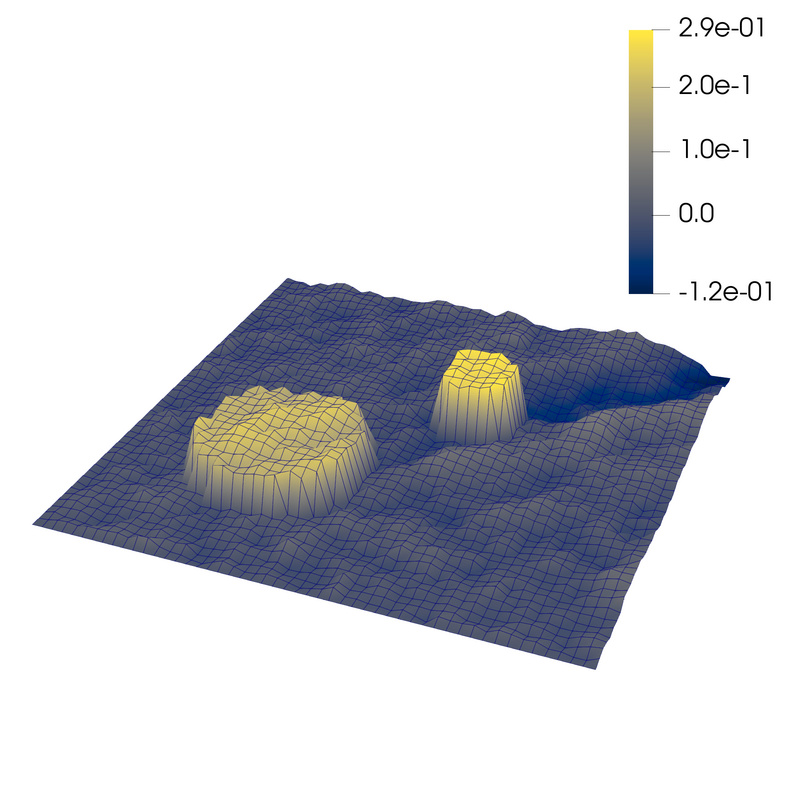}
    \caption{High-order reconstruction with optimal control for noisy data $\left( \sigma = 1\% \right)$.}
    \label{fig:2d_high_order_noise_oc_0.01}
\end{figure}

\begin{figure}
    \centering
    \begin{subfigure}{0.45\textwidth}
        \centering
        \includegraphics[width=\linewidth,trim={0 4cm 0 0.4cm},clip]{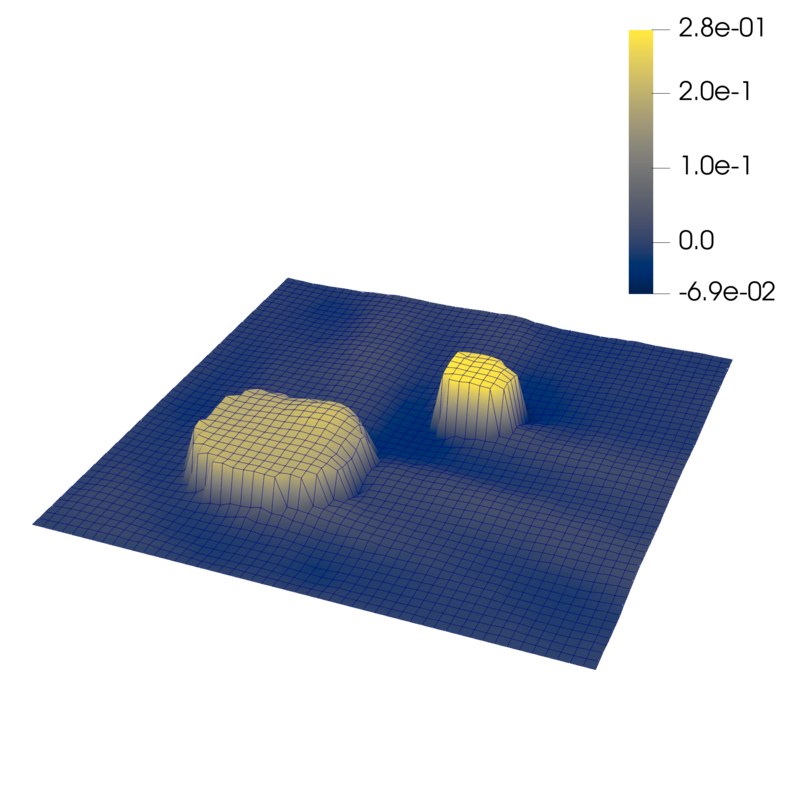}
        \caption{High-order reconstruction with TVD regularization for noisy data $\left( \sigma = 1\% \right)$}
        \label{fig:2d_high_order_noise_l2_0.01}
    \end{subfigure}
    \hfill
    \begin{subfigure}{0.45\textwidth}
        \centering
        \includegraphics[width=\linewidth,trim={0 4cm 0 0.4cm},clip]{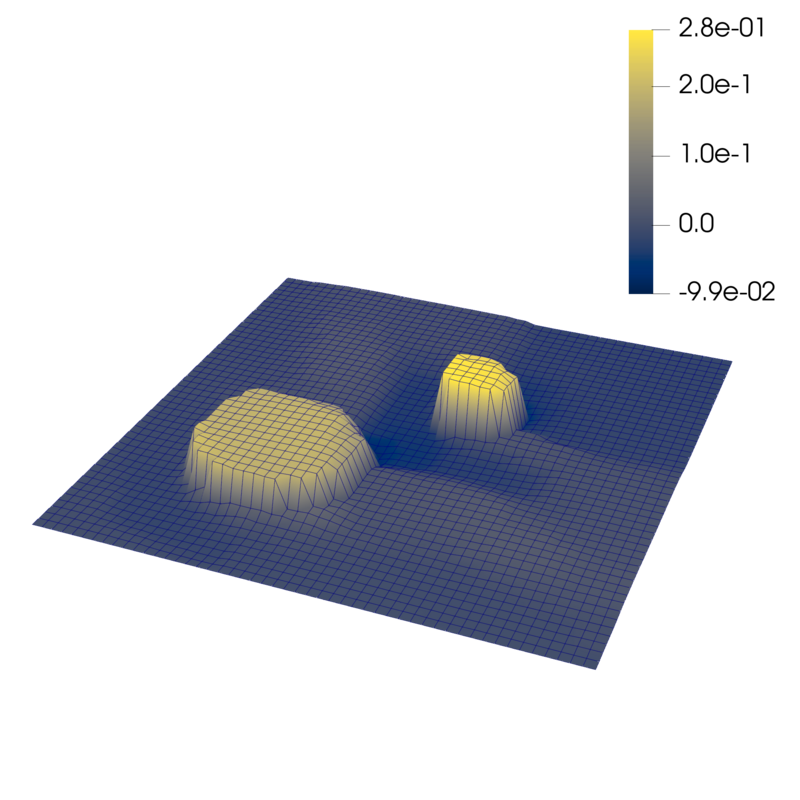}
        \caption{High-order reconstruction with $\Lonelone$ regularization for noisy data $\left( \sigma = 1\% \right)$}
        \label{fig:2d_high_order_noise_l1_0.01}
    \end{subfigure}
    \caption{High-order two-dimensional bathymetry reconstruction under $1\%$ Gaussian noise using TVD and $\Lonelone$ regularization.}
    \label{fig:2d_high_order_noise_0.01}
\end{figure}

We now examine the two-dimensional benchmark using the modified MCL solver at several noise levels. \Cref{fig:2d_high_order_noise_oc_0.01} shows that an unregularized optimal control reconstruction cannot cope with a $1\%$ noise level. It exhibits widespread oscillations throughout the domain. Applying a TVD penalty (\Cref{fig:2d_high_order_noise_l2_0.01}) removes these artifacts, yielding a smooth and more faithful estimate of the bathymetry. Using an $L^1(l_1)$ penalty (see \Cref{fig:2d_high_order_noise_l1_0.01}) yields a reconstruction similar to the TVD-stabilized solution. However, compared with the TVD case, the regions around the cylinder are noticeably flatter.

\begin{figure}
    \centering
    \begin{subfigure}{0.45\textwidth}
        \centering
        \includegraphics[width=\linewidth,trim={0 4cm 0 0.4cm},clip]{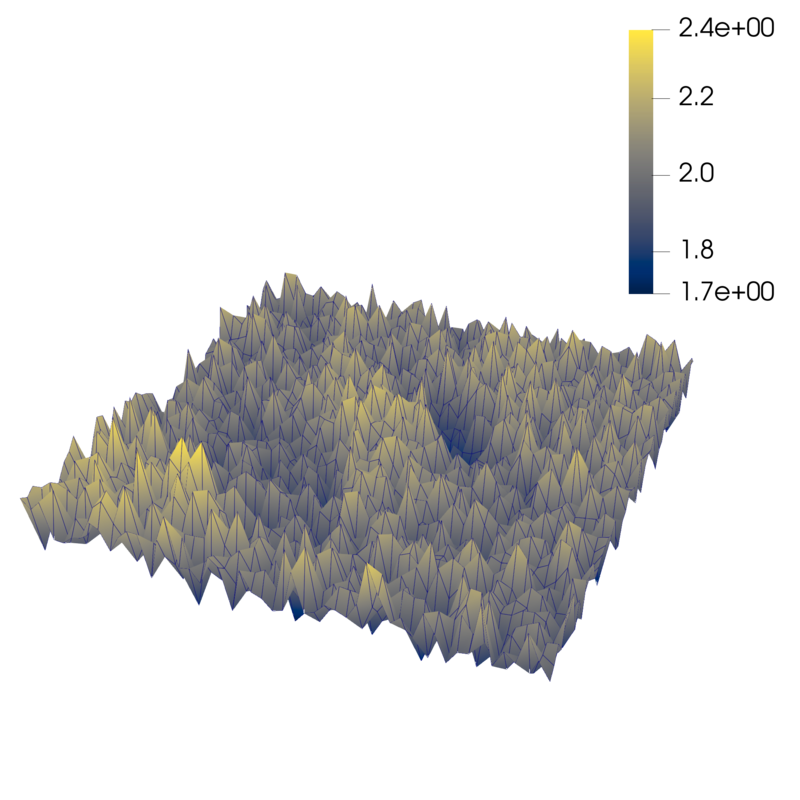}
        \caption{Reference free surface elevation at final time $\SI{60}{\second}$}
        \label{fig:2d_high_order_noise_total_height_0.05}
    \end{subfigure}
    \hfill
    \begin{subfigure}{0.45\textwidth}
        \centering
        \includegraphics[width=\linewidth,trim={0 4cm 0 0.4cm},clip]{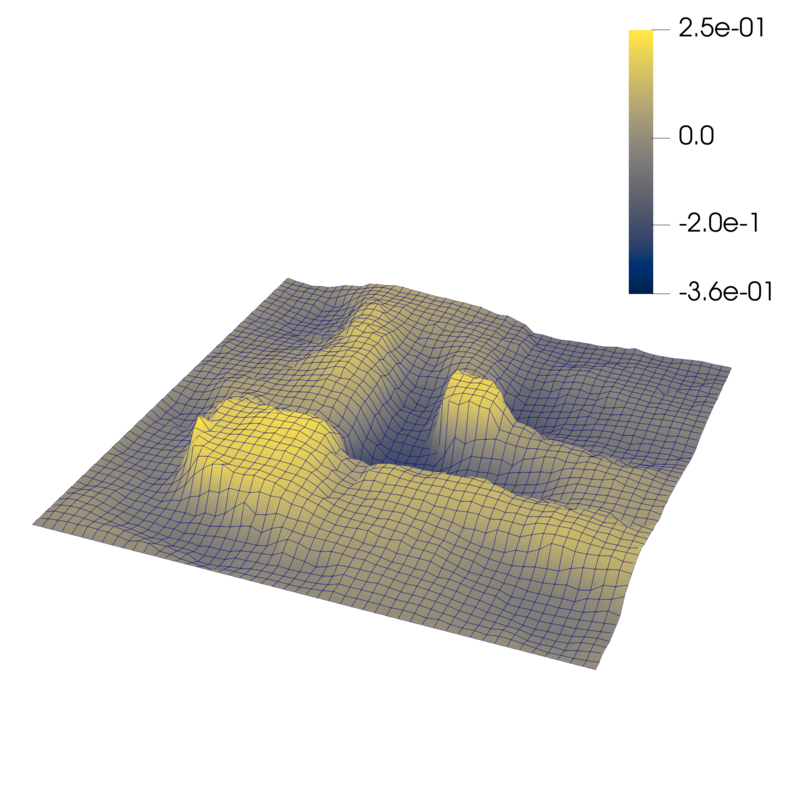}
        \caption{MCL reconstruction with TVD regularization}
        \label{fig:2d_high_order_noise_tvd_0.05}
    \end{subfigure}
    \begin{subfigure}{0.45\textwidth}
        \centering
        \includegraphics[width=\linewidth,trim={0 4cm 0 0.4cm},clip]{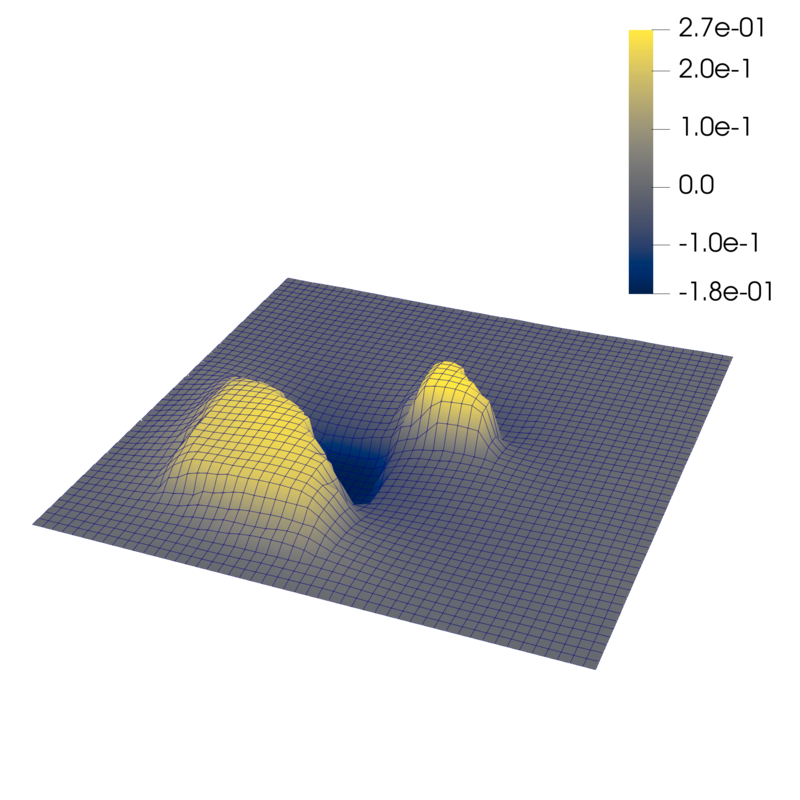}
        \caption{MCL reconstruction with $\Lonelone$ regularization}
        \label{fig:2d_high_order_noise_l1l1_0.05}
    \end{subfigure}
    \hfill
    \begin{subfigure}{0.45\textwidth}
        \centering
        \includegraphics[width=\linewidth,trim={0 4cm 0 0.4cm},clip]{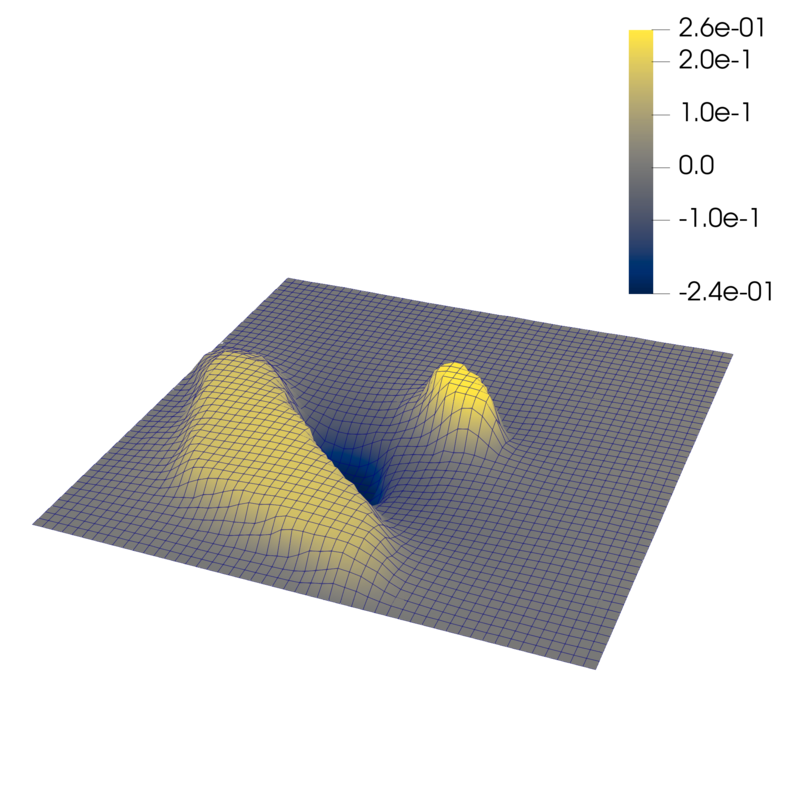}
        \caption{MCL reconstruction with $\Loneltwo$ regularization}
        \label{fig:2d_high_order_noise_l1l2_0.05}
    \end{subfigure}
    \caption{High-order two-dimensional bathymetry reconstructions for the two-cylinder benchmark under $5\%$ noise using TVD, $\Lonelone$, and $\Loneltwo$ regularization.}
    \label{fig:2d_high_order_noise_0.05}
\end{figure}

The more demanding test case corresponds to a $5\%$ additive noise level. \Cref{fig:2d_high_order_noise_total_height_0.05} illustrates a representative noisy measurement. \Cref{fig:2d_high_order_noise_tvd_0.05} displays the high-order reconstruction obtained with the TVD penalty term. Although the TVD regularizer produces a globally smooth surface, it introduces spurious humps that distort the true bathymetry. The maximum height is lower than expected. In contrast, the reconstruction obtained with the $L^1(l_1)$ regularization, shown in \Cref{fig:2d_high_order_noise_l1l1_0.05}, exhibits no oscillatory artifacts. The recovered bathymetry preserves all major features. The $L^1(l_2)$ variant in \Cref{fig:2d_high_order_noise_l1l2_0.05} behaves similarly to the $\Lonelone$ approach. However, the first cylinder is wider in orthogonal flow direction. Consequently, for this noise level the $L^1(l_1)$ method is preferable because it better preserves extreme values and avoids artificial humps. The overestimated gradient on the rear side of the first cylinder is reminiscent of that found in the one-dimensional result (see \Cref{fig:high_order_noise_0.01}).

\begin{table}
    \centering
    \begin{tabular}{cccccc}
       Algorithm ($\sigma$) & $L^2$ error (OC) & $L^2$ error \tvd & $L^2$ error (\Lonelone) & $L^2$ error (\Loneltwo) \\
       \midrule
       MCL (1\%) & \num{7.53e-1} & \num{5.85e-1} & \num{6.56e-1} & \num{6.70e-1} \\
       MCL (5\%) & \num{1.92e+0} & \num{1.98e+0} & \num{1.88e+0} & \num{2.05e+0} \\
       \bottomrule
    \end{tabular}
    \caption{Summary of $L^2$ errors for two-dimensional noisy data reconstructions using the MCL scheme with various regularizations.}
    \label{tab:noisy_error_2d}
\end{table}

For the two-dimensional noisy case (\Cref{tab:noisy_error_2d}), regularization is critical for stability. The optimal control baseline without gradient penalties yields errors exceeding \num{1.9} at 5\% noise, whereas TVD and $L^1$ variants stabilize the solution, reducing errors to approximately 0.6-0.7 at lower noise levels (1\%). However, under severe noise (5\%), all regularized methods exhibit increased error magnitudes relative to the clean case.

\FloatBarrier
\subsection{Large-scale simulations}\label{sec:large}
\begin{figure}
    \centering
    \begin{subfigure}{0.45\textwidth}
        \centering
        \includegraphics[width=\linewidth]{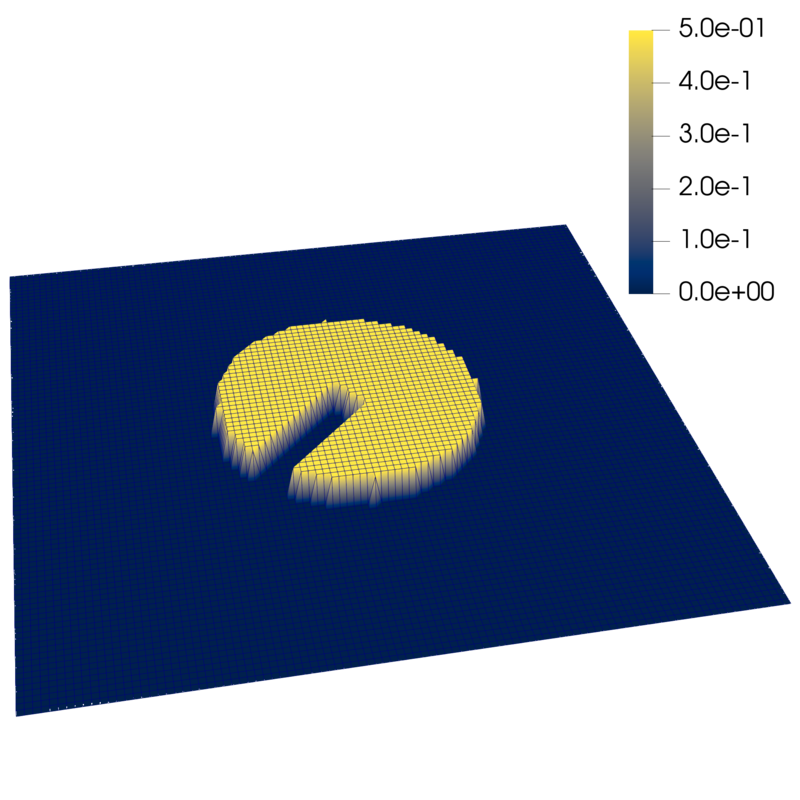}
        \caption{Interpolated initial condition of the sliced cylinder}
        \label{fig:2d_sliced_cylinder_init}
    \end{subfigure}
    \hfill
    \begin{subfigure}{0.45\textwidth}
        \centering
        \includegraphics[width=\linewidth]{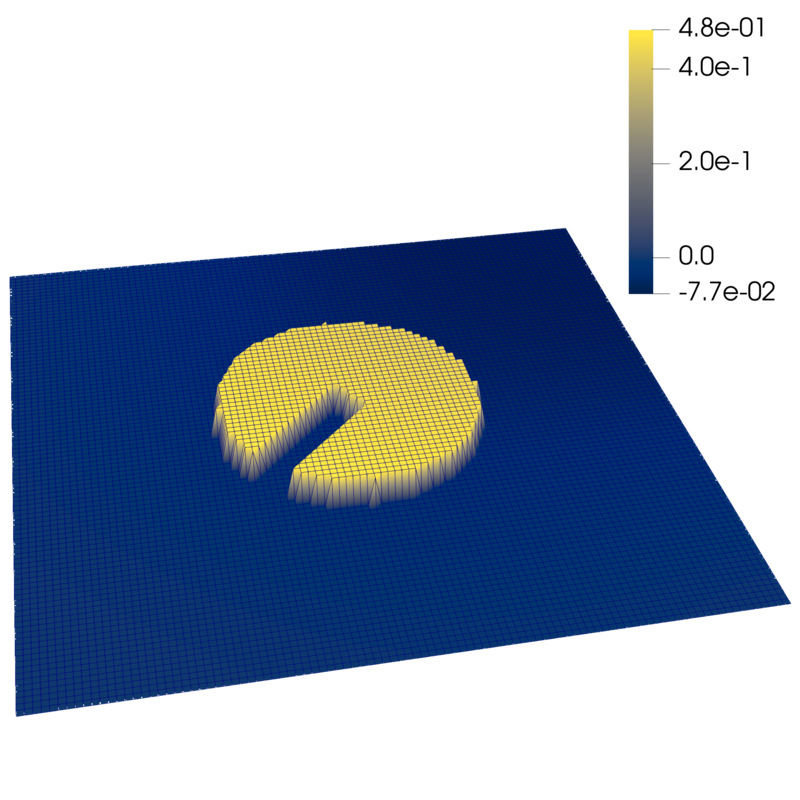}
        \caption{Low-order reconstruction with optimal control}
        \label{fig:2d_sliced_cylinder_oc}
    \end{subfigure}
    \caption{Large-scale sliced cylinder benchmark: Initial condition and low-order (ALF) optimal control reconstruction.}
    \label{fig:2d_sliced_cylinder}
\end{figure}

\noindent
To assess scalability on larger domains, we introduce an additional benchmark. The computational domain is $\Omega = [0,\SI{10}{\kilo\metre}]^2$, discretized by $100 \times 100$ quadrilateral elements ($\Delta x = \Delta y = \SI{100}{\metre}$). We prescribe a total height of $H = \SI{1000}{\metre}$ and a uniform velocity $\mathbf v = \left( 1.26, 1.26 \right)^{\top} \si{\metre\per\second}$, inspired by measurements presented in \citep{meinen2016structure}. The bathymetry is a sliced cylinder inspired by the solid-body rotation benchmark \citep{leveque1996high}. Let
\begin{align*}
    \mathbf x_s &= \mathbf x - 5000, \\
    f_{\text{circ}} (\mathbf x_s) &=
    \begin{cases}
         0.5, & \left| \mathbf x_s \right|_2 \leq 2000, \\
         0, &\text{otherwise,}
    \end{cases} \quad
    \mathbf{\tilde x} = \begin{pmatrix}
        \hphantom{+} \cos\left(\frac58 \pi\right) & \sin\left(\frac58 \pi\right) \\[.3em]
        -\sin\left(\frac58 \pi\right) & \cos\left(\frac58 \pi\right)
    \end{pmatrix} \mathbf x_s,\\
    f_{\text{slot}} (\mathbf{\tilde x}) &= \max \left\{ |\mathbf{\tilde x}_1| - 300, -\mathbf{\tilde x}_0 \right\}, \quad
    \tilde{b} = \begin{cases}
        f_{\text{circ}}, & f_{\text{slot}} < 0, \\
        0, &\text{otherwise}.
    \end{cases}
\end{align*}
The resulting bathymetry consists of a cylinder with a slot pointing upstream (see \Cref{fig:2d_sliced_cylinder_init}). The plots are scaled in the $z$-direction by \num{e3} for better visual representation. Using the optimal control framework with the ALF scheme, we recover the slice geometry. \Cref{fig:2d_sliced_cylinder_oc} shows that the reconstructed bathymetry matches the true bathymetry with a height of $\SI{0.5}{\metre}$. Sharp corners pose significant challenges due to large discrete gradients. Nevertheless, the objective functional \eqref{eq:J1} yields a stable solution for a fixed parameter of $\beta = \num{5e-5}$.

\begin{figure}
    \centering
    \begin{subfigure}{0.45\textwidth}
        \centering
        \includegraphics[width=\linewidth]{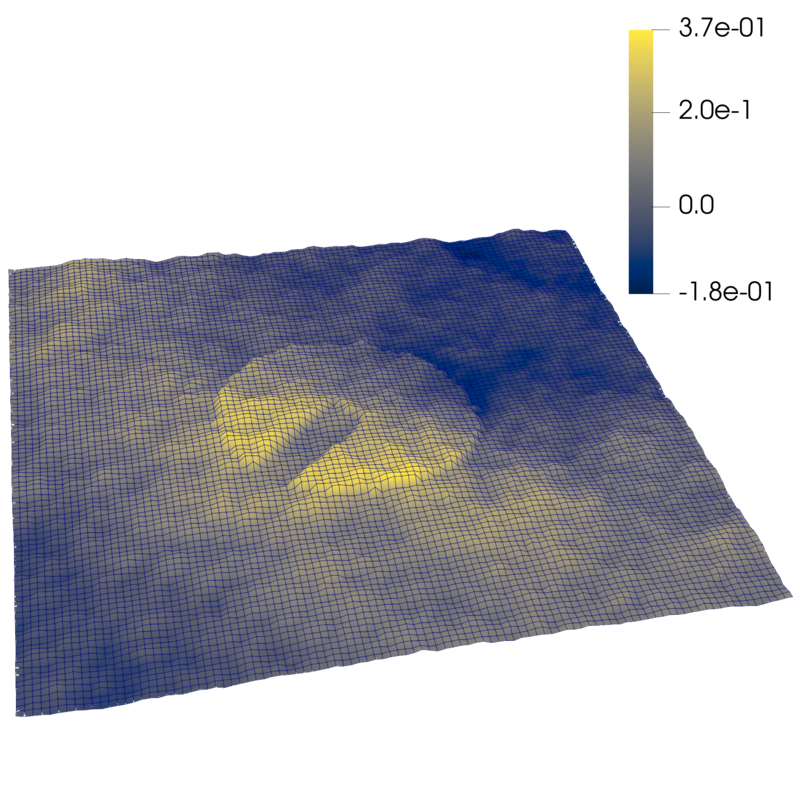}
        \caption{High-order reconstruction with $p=0$}
        \label{fig:2d_sliced_cylinder_with_noise_no}
    \end{subfigure}
    \hfill
    \begin{subfigure}{0.45\textwidth}
        \centering
        \includegraphics[width=\linewidth]{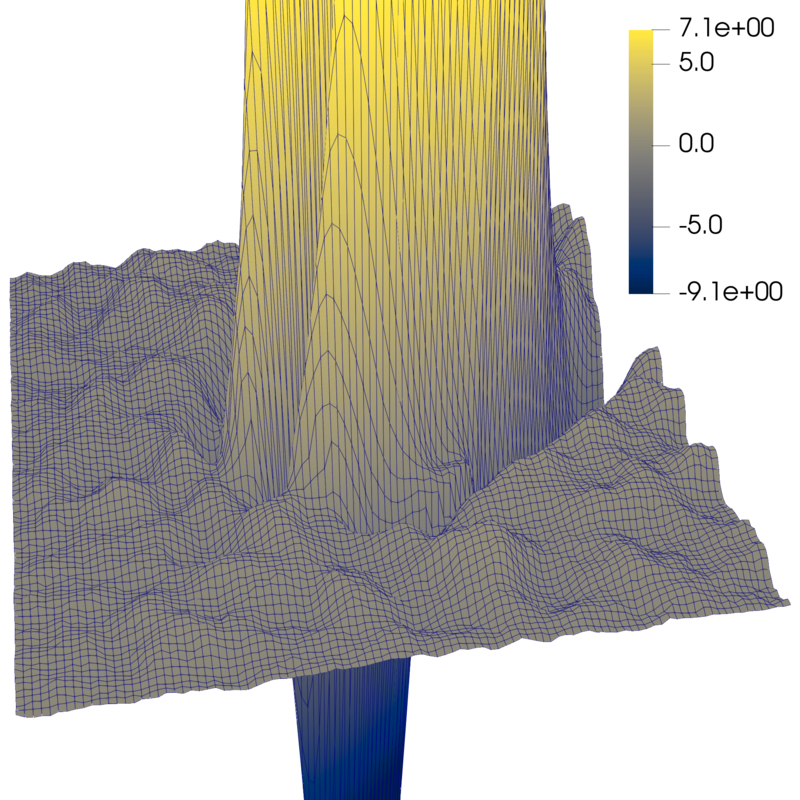}
        \caption{High-order reconstruction with OC}
        \label{fig:2d_sliced_cylinder_with_noise_oc}
    \end{subfigure}
    \begin{subfigure}{0.45\textwidth}
        \centering
        \includegraphics[width=\linewidth]{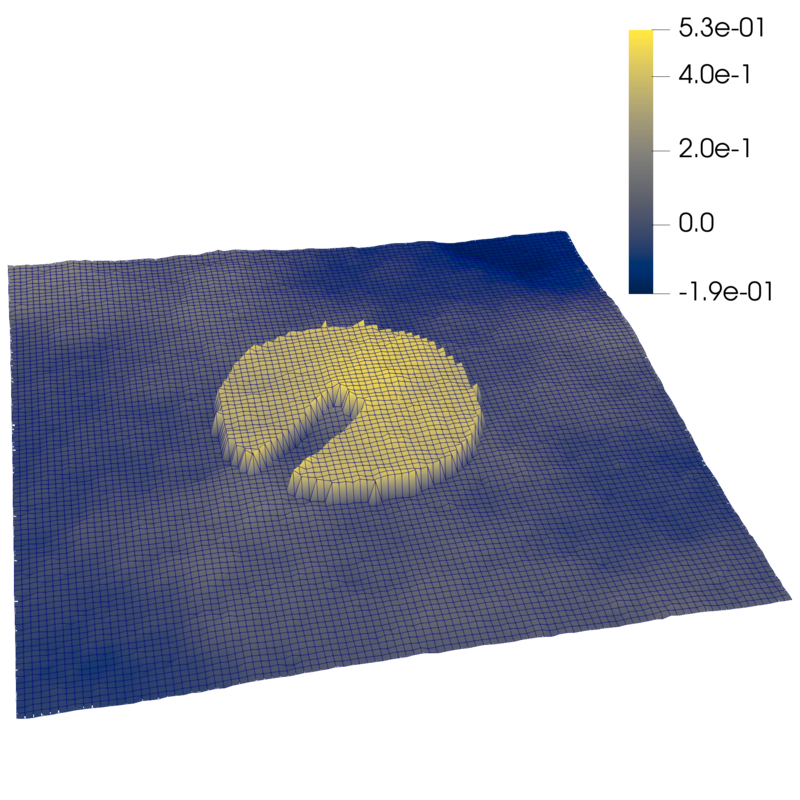}
        \caption{High-order reconstruction with $\Lonelone$ stabilization}
        \label{fig:2d_sliced_cylinder_with_noise_l1l1}
    \end{subfigure}
    \hfill
    \begin{subfigure}{0.45\textwidth}
        \centering
        \includegraphics[width=\linewidth]{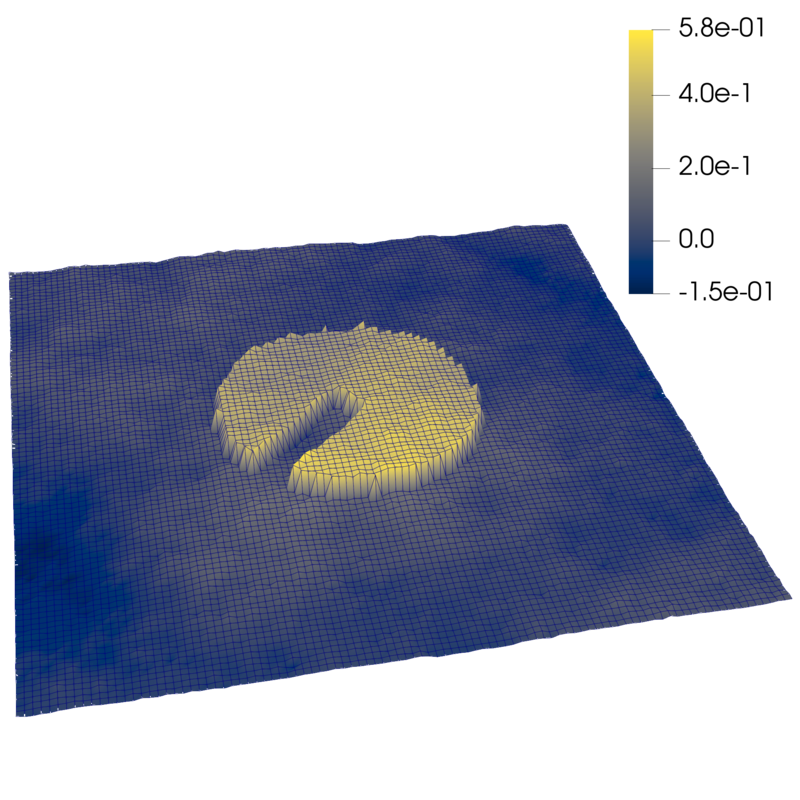}
        \caption{High-order reconstruction with $\Loneltwo$ regularization}
        \label{fig:2d_sliced_cylinder_with_noise_l1l2}
    \end{subfigure}
    \caption{Large-scale sliced cylinder benchmark: High-order reconstructions under noise using various stabilization techniques.}
    \label{fig:2d_sliced_cylinder_with_noise_reg}
\end{figure}

When noise with a standard deviation of $\num{3e-5}$ is added to the data and MCL is applied, the unstabilized reconstruction ($p=0$) deteriorates markedly \Cref{fig:2d_sliced_cylinder_with_noise_no}. The domain is covered with spurious oscillations and the maximum height is not recovered. \Cref{fig:2d_sliced_cylinder_with_noise_oc} demonstrates that, without a suitable gradient regularization, the optimal control solution fails to capture the sharp corners. The reconstructed bathymetry diverges and creates a hump that is significantly higher than the expected bathymetry.

These results underline the necessity of using spatially adaptive penalty terms, such as $L^1(l_1)$, to preserve discontinuities while suppressing noise in large-scale settings. The steep gradient observed at the beginning of the slice constitutes a major source of instability. When measurement noise is present, the $L^1(l_1)$ and $L^1(l_2)$ variants produce comparable results (see \Cref{fig:2d_sliced_cylinder_with_noise_l1l1,fig:2d_sliced_cylinder_with_noise_l1l2}). The cylindrical shape of the bathymetry is recovered, although the reconstructed extrema are slightly inaccurate.

\begin{table}
    \centering
    \begin{tabular}{ccccc}
       Algorithm ($\sigma$) & $L^2$ error (NOS) & $L^2$ error (OC) & $L^2$ error (\Lonelone) & $L^2$ error (\Loneltwo) \\
       \midrule
       ALF (0\%) & \num{4.45e+2} & \num{4.45e+2} & \num{3.94e+2} & \num{3.87e+2} \\
       MCL (3\%) & \num{1.52e+3} & \num{3.68e+4} & \num{8.23e+2} & \num{7.21e+2} \\
       \bottomrule
    \end{tabular}
    \caption{Summary of $L^2$ errors for large-scale sliced-cylinder reconstruction under noisy conditions.}
    \label{tab:large_error}
\end{table}

The error metrics in \Cref{tab:large_error} reveal significant sensitivity to regularization in large-scale settings. While the unstabilized high-order scheme yields errors in the order of \num{1.52e+3}, the optimal control variant without gradient penalties diverges further (\num{3.68e+4}), indicating numerical instability without sufficient gradient stabilization. In contrast, $L^1$-based regularizations reduce the error to approximately $\num{8e+2}$, demonstrating that sparsity-promoting terms are essential for maintaining stability and accuracy when scaling the inverse problem to larger computational meshes.

\section{Conclusion}\noindent
This study was focused on techniques for reconstructing bathymetry from measurement data within a flux-potential framework that enables optimal control at each time step. The initial approach involved solving an optimality system for the flux potentials. Although effective in consistency tests, this method encounters substantial challenges with noisy data. To address these limitations, we integrated a trust-region method into the optimization framework to handle non-smooth $L^1$ regularization. This allowed us to penalize the bathymetry gradient via an $L^1$ term, promoting sparsity and yielding superior results on perturbed data. Solving the smooth dual problem within a reduced-space formulation facilitated rapid convergence.

Performance is highly sensitive to the selection of optimization parameters. A key advantage of the proposed scheme is its adaptability: It accommodates any forward problem solver capable of computing time derivatives, making it suitable for both explicit and implicit methods. Although the Lax-Friedrichs approximation yielded promising results, its accuracy relative to analytical solutions remains limited. More accurate approaches, such as MCL, demonstrate better performance but require careful tuning of optimization parameters.

A primary structural advantage of this optimization-based stabilization is its efficient handling of streaming real-time data. Evaluation requires only the preceding state and current measurements, whereas full space-time adjoint methods would necessitate a complete backward integration over the entire dataset. The remaining challenges concern performance on large-scale simulations. Adaptive optimization strategies require further refinement. While noise-free scenarios permit looser constraints, noisy data result in a trade-off between accuracy and numerical stability.

Removing the bathymetry-dependent diffusion from the continuity equation improves the convergence behavior and the quality of reconstructions in our experiments. However, preservation of lakes at rest is no longer guaranteed. Although the optimization procedure tends to recover a balanced solution in practice, we have no proof that this is possible in all scenarios. Future research will focus on establishing a theoretical guarantee that optimality conditions enforce the discrete well-balancedness. The methods employed in the present paper require local optimality at individual time steps. Extensions to full space-time adjoint formulations appear feasible. Introducing adaptive time stepping could accelerate convergence in low-wave-speed regimes while enhancing resolution in high-wave-speed regions. Further improvements may involve incorporating friction forces or sediment transport physics, and exploring elastic-net regularization that blends the sparsity of $L^1$ with the smoothness of $L^2$. Additional investigations of regularization techniques on anisotropic meshes appear promising.

The numerical experiments performed so far utilized coarse meshes. To enable large-scale simulations, future work will explore the potential of GPU-accelerated forward problem solvers. Preliminary results demonstrate a significant speedup $\approx 10\times$. To further reduce the computational cost of data-driven bathymetry reconstruction, we intend to incorporate inexact evaluation strategies into the optimization framework. In particular, we aim to improve the algorithmic efficiency without compromising the quality of the results by using relaxed tolerances for intermediate iterates.

\section*{Acknowledgments}\noindent
This research received funding from the German Research Foundation (DFG) under grant number KU 1530/29-1. The authors are grateful to Dr. Hennes Hajduk (TU Dortmund University) for his expert guidance on the shallow-water equations and for providing constructive feedback that enhanced this work. Furthermore, the first author expresses his sincere gratitude to Dr. Drew Kouri and Dr. Denis Ridzal (Sandia National Laboratories) for guidance on integrating sparse optimization techniques into the flux-potential framework.

\begin{appendices}
\section{Numerical parameters}\label{app:parameter}\noindent
All numerical experiments in the main text were carried out with the parameter sets listed below. The tables are grouped by experiment type so that the reader can reproduce the benchmark computations in \Cref{sec:conv,sec:two,sec:noisy,sec:large}. The one-dimensional simulations were run on a single thread, the two-dimensional tests used 4 MPI ranks and were also benchmarked on Nvidia GPUs (performance evaluations are deferred to future work).

\begin{table}[!ht]
    \centering
    \begin{tabular}{l|c|c}
        Parameter & ALF (low-order) & MCL (high-order) \\
        \midrule
        Domain size & \multicolumn{2}{c}{$[0,\SI{25}{m}]$}\\
        Spatial discretization & \multicolumn{2}{c}{Equidistant segments on a uniform mesh}\\
        Mesh sizes $\Delta x$ & \multicolumn{2}{c}{$\frac14,\;\frac18,\;\frac1{16},\;\frac1{32}$}\\
        Time step $\Delta t$ & \multicolumn{2}{c}{$4\,\Delta x\,\SI{3e-2}{s}$}\\
        Final time $T$ & \multicolumn{2}{c}{$\SI{200}{s}$}\\
        Regularization weights & $\alpha=1,\;\beta=\num{e-11},\;\gamma=\num{e5}$ & $\alpha=1,\;\beta=\num{e-4},\;\gamma=\num{e5}$\\
        Solver tolerance (rel.) & \multicolumn{2}{c}{Intel's Cluster Sparse Solver}\\
        \bottomrule
    \end{tabular}
    \caption{Numerical parameters for the one-dimensional convergence benchmark (\Cref{sec:conv}).}
    \label{tab:conv}
\end{table}

\begin{table}[!ht]
    \centering
    \begin{tabular}{l|c}
        Parameter & MCL (high-order)\\
        \midrule
        Domain size & $[0,\SI{25}{m}]^2$\\
        Mesh & $50\times 50$ quadrilaterals ($\Delta x=\Delta y=\SI{0.5}{m}$)\\
        Time step $\Delta t$ & $\SI{e-2}{s}$\\
        Final time $T$ & $\SI{60}{s}$\\
        Regularization weights & $\alpha=1,\;\beta=\num{e-7},\;\gamma=1$\\
        Solver tolerance (rel.) & $\num{e-8}$ (MINRES + hypre's BoomerAMG)\\
        \bottomrule
    \end{tabular}
    \caption{Numerical parameters for the two-dimensional cylinder benchmark (\Cref{sec:two}).}
    \label{tab:two}
\end{table}

\begin{table}[!ht]
    \centering
    \begin{tabular}{l|c|c}
        Parameter & ALF (low-order) - 1\% noise & ALF (low-order) - 5\% noise\\
        \midrule
        Domain size & \multicolumn{2}{c}{$[0,\SI{25}{m}]$}\\
        Mesh & \multicolumn{2}{c}{$100$ segments ($\Delta x=\SI{0.25}{m}$)}\\
        Time step $\Delta t$ & \multicolumn{2}{c}{$\SI{3e-2}{s}$}\\
        Final time $T$ & \multicolumn{2}{c}{$\SI{200}{s}$}\\
        Regularization weights & \multicolumn{2}{c}{$\alpha=1,\;\gamma=\num{e5}$}\\
        Optimal control & \multicolumn{2}{c}{$\beta=\num{e-4}$}\\
        TVD stabilization & $\beta=\num{e-9},\;\epsilon=\num{e-2},\;\zeta=\num{e-4}$ & $\beta=\num{e-9},\;\epsilon=5\times10^{-2},\;\zeta=\num{e-4}$\\
        $L^1$ penalty & $\beta=\num{e-9},\;\kappa=\num{e-2},\;\nu=1$ & $\beta=\num{e-9},\;\kappa=5\times10^{-2},\;\nu=1$\\
        Solver tolerance (rel.) & \multicolumn{2}{c}{Intel's Cluster Sparse Solver}\\
        \bottomrule
    \end{tabular}
    \caption{Numerical parameters for one-dimensional noisy data experiments using the ALF scheme (\Cref{sec:noisy}).}
    \label{tab:noisy_1D_ALF}
\end{table}

\begin{table}[!ht]
    \centering
    \begin{tabular}{l|c}
        Parameter & MCL (high-order) - 1\% noise\\
        \midrule
        Domain size & $[0,\SI{25}{m}]$\\
        Mesh & $100$ segments ($\Delta x=\SI{0.25}{m}$)\\
        Time step $\Delta t$ & $\SI{3e-2}{s}$\\
        Final time $T$ & $\SI{200}{s}$\\
        Regularization weights & $\alpha=1,\;\gamma=\num{e5}$\\
        Optimal control & $\beta=\num{e-4}$\\
        TVD stabilization & $\beta=\num{e-9},\;\epsilon=\num{e-2},\;\zeta=\num{e-4}$\\
        $L^1$ penalty & $\beta=\num{e-9},\;\kappa=0.015,\;\nu=\num{e-1}$\\
        Solver tolerance (rel.) & Intel's Cluster Sparse Solver\\
        \bottomrule
    \end{tabular}
    \caption{Numerical parameters for one-dimensional noisy data experiments using the MCL scheme (\Cref{sec:noisy}).}
    \label{tab:noisy_1D_MCL}
\end{table}

\begin{table}[!ht]
    \centering
    \begin{tabular}{l|c|c}
        Parameter & MCL (high-order) - 1\% noise & MCL (high-order) - 5\% noise\\
        \midrule
        Domain size & \multicolumn{2}{c}{$[0,\SI{25}{m}]^2$}\\
        Mesh & \multicolumn{2}{c}{$50\times 50$ quadrilaterals ($\Delta x=\Delta y=\SI{0.5}{m}$)}\\
        Time step $\Delta t$ & \multicolumn{2}{c}{$\SI{e-2}{s}$}\\
        Final time $T$ & \multicolumn{2}{c}{$\SI{60}{s}$}\\
        Regularization weights & \multicolumn{2}{c}{$\alpha=1,\;\gamma=\num{e5}$}\\
        Optimal control & \multicolumn{2}{c}{$\beta=\num{e-4}$}\\
        TVD stabilization & $\beta=\num{e-9},\;\epsilon=\num{e-2},\;\zeta=\num{e-4}$ & $\beta=\num{e-9},\;\epsilon=5\times10^{-2},\;\zeta=5\times10^{-4}$\\
        $L^1$ penalty & $\beta=\num{e-9},\;\kappa=\num{e-2},\;\nu=\num{e-1}$ & $\beta=\num{e-9},\;\kappa=0.15,\;\nu=\num{e-1}$\\
        Solver tolerance (rel.) & \multicolumn{2}{c}{$\num{e-8}$ (GMRES + hypre's Euclid (ILU))}\\
        \bottomrule
    \end{tabular}
    \caption{Numerical parameters for two-dimensional noisy data experiments using the MCL scheme (\Cref{sec:noisy}).}
    \label{tab:noisy_2D}
\end{table}

\begin{table}[!ht]
    \centering
    \begin{tabular}{l|c}
        Parameter & ALF (low-order) \& MCL (high-order)\\
        \midrule
        Domain size & $[0,\SI{10}{km}]^2$\\
        Mesh & $100\times 100$ quadrilaterals ($\Delta x=\Delta y=\SI{100}{m}$)\\
        Time step $\Delta t$ & $\SI{e-1}{s}$\\
        Final time $T$ & $\SI{1}{h}$\\
        Regularization weights & $\alpha=1,\;\beta=5\times10^{-5},\;\gamma=1$\\
        $L^1$ penalty & $\kappa=1,\;\nu=5\times10^{-5}$\\
        Solver tolerance (rel.) & $\num{e-8}$ (MINRES + hypre's BoomerAMG)\\
        \bottomrule
    \end{tabular}
    \caption{Numerical parameters for the large-scale sliced cylinder benchmark (\Cref{sec:large}).}
    \label{tab:large}
\end{table}
\end{appendices}

\FloatBarrier
\bibliographystyle{elsarticle-num}
\bibliography{literature}

\end{document}